\documentclass{amsart}
\usepackage{tikz}
\usepackage{xcolor}
\usepackage{amssymb,latexsym,amsmath,extarrows}
\usepackage{graphicx,mathrsfs}
\usepackage{hyperref}

\numberwithin{equation}{section}

\newtheorem{theorem}{Theorem}[section]
\newtheorem{lemma}[theorem]{Lemma}

\newtheorem{proposition}[theorem]{Proposition}
\newtheorem{remark}[theorem]{Remark}

\newtheorem{definition}[theorem]{Definition}
\newtheorem{corollary}[theorem]{Corollary}

\renewcommand{\Bbb}{\mathbb}

\newcommand{\al}{\alpha}
\newcommand{\be}{\beta}

\newcommand{\Ga}{\Gamma}
\newcommand{\de}{\delta}

\newcommand{\De}{\Delta}
\newcommand{\e}{\varepsilon}

\newcommand{\ka}{\kappa}
\newcommand{\la}{\lambda}

\newcommand{\La}{\Lambda}
\newcommand{\si}{\sigma}
\newcommand{\Si}{\Sigma}
\newcommand{\vp}{\varphi}
\newcommand{\om}{\omega}

\newcommand{\cp}{\mathcal P}

\newcommand{\cn}{\mathcal N}

\newcommand{\ce}{\mathcal E}

\newcommand{\wt}{\widetilde}
\newcommand{\wh}{\widehat}
\newcommand{\lp}{L^{p}}

\newcommand{\ZR}{\mathbb{R}}
\newcommand{\ZT}{\mathbb{T}}
\newcommand{\ZZ}{\mathbb{Z}}
\newcommand{\ZP}{\mathbb{P}}
\newcommand{\ZC}{\mathbb{C}}

\newcommand{\ZN}{\mathbb{N}}
\newcommand{\ZS}{\mathbb{S}}

\newcommand{\Id}{{\bf 1}}

\newcommand{\cl}{{\mathcal L}}

\newcommand{\rapid}{\mathrm{RapDec}}
\newcommand{\dist}{\mathrm{dist}}
\newcommand{\supp}{\mathrm{supp}}

\begin{document}

\title[Refined Strichartz and maximal extension operator]{A note on the refined Strichartz estimates and maximal extension operator}
\author{Shukun Wu}

\address{
Shukun Wu\\
Department of Mathematics\\
University of Illinois at Urbana-Champaign\\
Urbana, IL, 61801, USA}

\email{shukunw2@illinois.edu}
\date{}

\begin{abstract}
There are two parts for this paper. In the first part we extend some results in a recent paper by Du, Guth, Li and Zhang to a more general class of phase functions. The main methods are Bourgain-Demeter's $l^2$ decoupling theorem and induction on scales. In the second part we prove some positive results for the maximal extension operator for hypersurfaces with positive principal curvatures. The main methods are sharp $L^2$ estimates by Du and Zhang, and the bilinear method by Wolff and Tao.
\end{abstract}
\maketitle

\section{Introduction}
\setcounter{equation}0
Let $\Phi(\xi):B^n(0,2)\to\ZR$ be a real smooth function with $|\nabla\Phi|\leq10$ and $\det(\nabla^2\Phi)\not=0$. Typical examples for $\Phi(\xi)$ are: $\Phi(\xi)=\sqrt{9-|\xi|^2}$, the truncated sphere; $\Phi(\xi)=v_1\xi_1^2+\cdots+v_n\xi_n^2$, $v_j\in\{-1,1\}$, the truncated hyperbolic paraboloid. Let $f:\ZR^n\to\ZC$ be such that ${\rm supp}(\wh{f})\subset B^n(0,1)$. We consider the extension operator $\ce_\Phi$ associated to the hypersurface $(\xi,\Phi(\xi))$ that is defined as
\begin{equation}
\label{extension}
    \ce_{\Phi} f(x,t):=\int_{\ZR^n}e^{i(x\cdot\xi+t\Phi(\xi))}\wh{f}(\xi)d\xi,
\end{equation}
and the maximal operator $(\sup_{t\in\ZR}\ce_{\Phi})$ defined as
\begin{equation}
    \label{maximal}
    (\sup_{t\in\ZR}\ce_{\Phi})f=\sup_{t\in\ZR}\Big|\int_{\ZR^n}e^{i(x\cdot\xi+t\Phi(\xi))}\wh{f}(\xi)d\xi\Big|.
\end{equation}

The operator $\ce_\Phi$ can be viewed as a generaliztion of Schr\"odinger operator: $\ce_{\Phi}$ coincides with the Schr\"odinger operator $e^{it\De}$ when $\Phi(\xi)=|\xi|^2$. In particular, $e^{it\De}f$ is the solution of the free Schr\"odinger equation
\begin{equation}
    \left\{
    \begin{array}{ll}
    iu_t-\Delta u=0, & (x,t)\in\ZR^n\times\ZR^+.  \\[1ex]
    u(x,0)=f(x),  &  x\in\ZR^n.
    \end{array}
    \right.
\end{equation}

There is a long history for the study of Maximal Schr\"odinger operator $\sup_{t>0}e^{it\De}$ and the related topics. In 1979, Carleson \cite{Carleson-Schrodinger} considered the following problem: Identify the range of $s$ such that $\lim_{t\to0^+}e^{it\De}f(x)=(2\pi)^nf(x)$ almost everywhere whenever $f\in H^s(\ZR^n)$. In the same paper, Carleson proved $s\geq1/4$ when $n=1$, and Dahlberg and Kenig showed the result is sharp in \cite{Dahlberg-Kenig}.

When $n\geq2$, Bourgain gave counterexamples in \cite{Bourgain-counterexample} indicating that the almost everywhere convergence problem can fail if $s<\frac{n}{2(n+1)}$. Several months later, Du, Guth and Li \cite{Du-Li-Guth-Schrodinger-R3} proved that $\|\sup_{0<t\leq1}|e^{it\De} f|\|_3\leq\|f\|_{H^s}$ when $n=2$ and $s>\frac{1}{3}$. These $(H^s,L^3)$ estimates automatically imply the sharp bound for the convergence problem in $\ZR^2$ up to endpoints. Later, using multilinear refined Strichartz estimates that was first introduced in \cite{Du-Li-Guth-Schrodinger-R3}, Du, Guth, Li and Zhang \cite{Du-Li-Guth-Zhang-Schrodinger-Rn} prove that the almost convergence problem holds for $s>\frac{n+1}{2(n+2)}$ when $n\geq3$. Recently, up to endpoints, the convergence problem was completely settled by Du and Zhang \cite{Du-Zhang} in higher dimensions. Their method is partially based on the multilinear refined Strichartz estimates.

\medskip

In the first part of the paper, we extend the results in \cite{Du-Li-Guth-Zhang-Schrodinger-Rn} to the extension operator $\ce_\Phi$ defined in \eqref{extension}. Before stating our results, we make a definition concerning the principal curvatures of a smooth hypersurface.
\begin{definition}
Let $\Si$ be a smooth hypersurface in $\ZR^{n+1}$ with nonzero Gaussian curvature. We denote by $d(\Si)$ the minimum of the number of positive and negative principal curvatures. We let
$d(\Phi)=d(\Si)$ if $\Si=(\xi,\Phi(\xi))$ is the graph $\Phi$.
\end{definition}
\begin{theorem}
\label{refined-Strichartz}
Let $\ce_\Phi$ be defined in \eqref{extension}. Let   $\frac{2(n+2-d(\Phi))}{n-d(\Phi)}\leq p<\infty$ and let ${\bf B}=\{B_j\}$ be a collection of finitely overlapping $R^{1/2}$-cubes in $B^{n+1}(0,R)$, such that
\begin{equation}
    \label{spreading-out-condition}
    \|\ce_{\Phi} f\|_{L^p(B_j)}\sim\|\ce_{\Phi} f\|_{L^p(B_{j'})}.
\end{equation}
We further assume that $\bf B$ has the "spreading out" property. That is, those $B_j\in {\bf B}$ are arranged in horizontal slabs of the form $\ZR^n\times\{t_0,t_0+R^{1/2}\}$, so that each slab contains $\sim\sigma$ cubes $B_j$. Define $Y=\cup_jB_j$. Then, for any $\e>0$, there is  a constant $C_\e$ such that
\begin{equation}
\label{refined-Strichartz-esti}
    \|\ce_{\Phi} f\|_{L^p(Y)}\leq C_{p,\e} R^\e\si^{-\frac{p-2}{2p}}\|f\|_{L^2}.
\end{equation}
\end{theorem}

\begin{remark}
\rm

The endpoint $\frac{2(n+2-d(\Phi))}{n-d(\Phi)}$ in Theorem \ref{refined-Strichartz} is the best one can get. For example, let us take $\Phi:\ZR^2\to\ZR,~\Phi(\xi)=\xi_1\xi_2$, so $d(\Phi)=1$. Consider a test function $\wh{f}(\xi):=\vp(R\xi_1)\vp(\xi_2)$, where $\vp:[-2,2]\to\ZR$ is a bump function of the interval $[-1,1]$, so $\|f\|_2\sim R^{-1/2}$. This implies
\begin{equation}
    \ce_{\Phi} f(x_1,x_2,t)=R^{-1}\int_{\ZR^2}e^{iR^{-1}x_1\xi_1}e^{ix_2\xi_2}e^{iR^{-1}t\xi_1\xi_2}\vp(\xi_1)\vp(\xi_2)d\xi_1d\xi_2.
\end{equation}
Since the gradient of the phase function is $(R^{-1}x_1+R^{-1}t\xi_2,~x_2+R^{-1}t\xi_1)$, one can use the method of (non) stationary phase to conclude that
$\|\ce_\Phi f\|_{L^p(B)}\sim R^{1/p-1}$ for any $R^{1/2}$ cube $B$ in the thin slab $Y:=\{(x_1,x_2,t):|x_1|,|t|\leq R/2, |x_2|\leq R^{1/2}\}$. Hence $\|\ce_\Phi f\|_{L^p(Y)}\sim R^{2/p-1}$. But in the set $Y$, each horizontal slab $\{(x,t):|t-t_0|\leq R^{1/2},|t_0|\leq R/2\}$ contains $\sim R^{1/2}$ many $R^{1/2}$ cubes. So if \eqref{refined-Strichartz-esti} was true, one needs
\begin{equation}
    R^{2/p-1}\sim\|\ce_\Phi f\|_{L^p(Y)}\leq C_\e R^\e R^{-\frac{p-2}{4p}}\|f\|_2\approx R^{-\frac{p-2}{4p}}R^{-1/2},
\end{equation}
which is true only when $p\geq 6=\frac{2(n+2-d(\Phi))}{n-d(\Phi)}$. We can construct similar examples in higher dimensions to show that the endpoint $\frac{2(n+2-d(\Phi))}{n-d(\Phi)}$ is sharp.

\end{remark}

Similar to the multilinear results in \cite{Du-Li-Guth-Zhang-Schrodinger-Rn}, we have a multilinear analogue of Theorem \ref{refined-Strichartz}. Before stating the theorem, we need the following definition about the transversity of a collection of functions.

\begin{definition}
Let $f_j:\ZR^n\to\ZC$, $j=1,2,\ldots,k$. We say $\{f_j\}$ have frequencies $k$-transversely supported in $B^n(0,1)$, if for any $\xi\in{\rm supp}(\wh{f})\subset B^n_1$,
\begin{equation}
\label{k-transverse}
    |V_A(\xi_1)\wedge\cdots\wedge V_A(\xi_k)|\geq c>0,
\end{equation}
where $c$ is an absolute constant, and $V_A(\xi)=(\nabla\Phi(\xi),-1)$.
\end{definition}

With the definition of "frequencies $k$-transversely supported" in hand, we now can state our multilinear results.

\begin{theorem}
\label{multilinear-refined-strichartz}
Let $\ce_\Phi$ be defined in \eqref{extension} and $\frac{2(n+2-d(\Phi))}{n-d(\Phi)}\leq p<\infty$ and $2\leq k\leq n+1$. Suppose $f_j:\ZR^n\to\ZC$ have frequencies $k$-transversely supported in $B^n(0,1)$. Suppose $B_1,B_2,\ldots,B_N$ are finitely overlapping $R^{1/2}$-cubes in $B^{n+1}(0,R)$, such that for each $B_l$,
\begin{equation}
    \|\ce_{\Phi} f_j\|_{L^p(B_l)}\sim\|\ce_{\Phi} f_j\|_{L^p(B_{l'})}
\end{equation}
for all $j=1,2,\ldots,k$. Let $Y=\cup_{l=1}^N B_l$. Then for any $\e>0$, there exists an absolute constant $C_\e$ so that
\begin{equation}
    \label{multilienar-refined-strichartz-esti}
    \Big\|\prod_{j=1}^k|\ce_{\Phi}f_j|^{1/k}\Big\|_{\lp(Y)}\leq C_\e R^\e N^{-\frac{(k-1)(p-2)}{2kp}}\prod_{j=1}^k\|f_j\|_2^{1/k}.
\end{equation}
\end{theorem}

\medskip

In the second part of the paper, we study the $(H^s,L^p)$ behavior for the maximal extension operator \eqref{maximal} in $n\geq3$, with $d(\Phi)=0$. Here we require $\Phi:\ZR^n\to\ZR$ to be a global function. Consider the following problem: Let $f:\ZR^n\to\ZC$ and $\e>0$. For $p>2$, find the optimal $p$ and $s$ such that
\begin{equation}
\label{Maximal-schrodinger-conj}
    \big\|\sup_{t\in\ZR }|\ce_{\Phi}f|\big\|_{L^p(\ZR^n)}\leq C_{p,\e} \|f\|_{H^{s+\e}}. 
\end{equation}
Note that when $\wh{f}$ is supported in the annulus $\{|\xi|\sim R\}$ for $R\geq1$, the estimate \eqref{Maximal-schrodinger-conj} boils down to  
\begin{equation}
\label{scale-R}
    \big\|\sup_{t\in\ZR }|\ce_{\Phi}f|\big\|_{L^p(\ZR^n)}\leq C_{p,\e} R^{n/2-n/p}\|f\|_2.
\end{equation}
This shows that the optimal $s$ one could expect is $s=n/2-n/p$. Let assume $s=n/2-n/p$ from now on.

There is one more condition we need to impose on $\Phi$. Note that when recaling the annulus $|\xi|\sim R$ back to $|\xi|\sim 1$ in \eqref{scale-R}, we have to deal with the scaled function $\Phi_R(\xi):=R^{-2}\Phi(R\xi)$. Hence, we require $\Phi$ to satisfy that for any $\al\in\ZN^n$ and $|\xi|\geq1$, 
\begin{equation}
\label{derivative-phi}
    |\partial^\al\Phi(\xi)|\leq C_{\al}|\xi|^{2-|\al|}.
\end{equation}
In fact, as shown in \cite{Hormander} Theorem 7.7.5, we only need \eqref{derivative-phi} for $|\al|\leq Cn$ so that one could obtain the decay estimate for $\ce_{\Phi_R}\phi$ uniformly in $R$, where $\phi$ is a smooth function whose Fourier transform is supported in $B^n(0,2)$. We will use the estimate for $\ce_{\Phi}\phi$ in Section 6.


Due to Bourgain's counterexample in \cite{Bourgain-counterexample}, one might conjectures that \eqref{Maximal-schrodinger-conj} is true for $\Phi(\xi)=|\xi|^2$, $p>\frac{2(n+1)}{n}$ and $n\geq3$. However, recently Du, Kim, Wang and Zhang gave new counterexamples in \cite{Du-Kim-Wang-Zhang} indicating that \eqref{Maximal-schrodinger-conj} can only holds for 
\begin{equation}
p\geq\max_{m\in\ZN,1\leq m\leq n} 2+\frac{4}{n-1+m+n/m}.
\end{equation}

We give some positive results to the estimate \eqref{Maximal-schrodinger-conj} via bilinear techniques introduced by Wolff \cite{Wolff} and Tao \cite{Tao-Bilinear}.

\begin{theorem}
\label{Sobolev}
Let $\Phi:\ZR^n\to\ZR$ be a function satisfying $d(\Phi)=0$, $|\nabla\Phi|\lesssim1$ and $|\det(\nabla^2\Phi)|\sim1$. Suppose that in addition $\Phi$ satisfies \eqref{derivative-phi}.  We let $\ce_\Phi f$ be defined in \eqref{extension}, and let $f:\ZR^n\to\ZC$. Then for any $p>2+\frac{4}{n+2-1/n}$,  $s=\frac{n}{2}-\frac{n}{p}$ and $\e>0$, there is a constant $C_\e$ only depends on dimension and $\e$, so that
\begin{equation}
\label{Sobolev-esti}
    \big\|\sup_{t\in\ZR }|\ce_{\Phi}f|\big\|_{L^p(\ZR^n)}\leq C_{p,\e} \|f\|_{H^{s+\e}}.
\end{equation}
\end{theorem}
\noindent One can partition $f$ on the Fourier side using Littlewood-Paley decomposition so that Theorem \ref{Sobolev} reduces to the following estimates:

\begin{theorem}
\label{Maximal-Schrodinger}
Suppose that $\Phi:B^n(0,2)\to\ZR$ is a function satisfying $d(\Phi)=0$, $|\nabla\Phi|\lesssim1$ and $|\det(\nabla^2\Phi)|\sim1$. Let $f:\ZR^n\to\ZC$ be with $\wh{f}\in L^2(B^n(0,1))$. Then for any $p>2+\frac{4}{n+2-1/n}$, there is a constant $C$ only depends dimension, so that
\begin{equation}
\label{Maximal-Schrodinger-esti}
    \big\|\sup_{t\in\ZR }|\ce_{\Phi}f|\big\|_{L^p(\ZR^n)}\leq C_p \|f\|_2.
\end{equation}
\end{theorem}

\noindent In section 4, we will see that the maximal extension operator $(\sup_{t}\ce_\Phi)$ is indeed a local version of the classic extension operator $\ce_\Phi$. Hence it is reasonable to believe \eqref{Maximal-Schrodinger-esti} can hold for $p$ beyond the range of classic Strichartz estimates.

Another interesting problem would be : Under the setting of Theorem \ref{Maximal-Schrodinger}, find the smallest $p$ such that
\begin{equation}
\label{Maximal-Extension-pp}
     \big\|\sup_{t\in\ZR }|\ce_{\Phi}f|\big\|_{L^p(\ZR^n)}\leq C_p \|\wh{f}\|_p.
\end{equation}
From the same reason mentioned above, the best $p$ for \eqref{Maximal-Extension-pp} is expected to be smaller than $\frac{2(n+1)}{n}+\e$, the endpoint for restriction conjecture. Notice that \eqref{Maximal-Schrodinger-esti} automatically implies $\eqref{Maximal-Extension-pp}$ for $p>2+\frac{4}{n+2-1/n}$. We will not attack \eqref{Maximal-Extension-pp} in this paper.

Recently in \cite{Du-Guth-Ou-Wang-Wilson-Zhang}, the authors there used refined Strichartz estimates to study weighted restriction problems. Instead of using bilinear techniques, they rely on the method of polynomial partitioning, which was introduced to the study of oscillatory integrals by Guth in \cite{Guth-restriction-R3}. It is not implausible that one can use their methods and techniques to obtain further improvements on Theorem \ref{Maximal-Schrodinger}, while we do not pursue it here.

\medskip

This paper is organized as follows: We prove Theorem \ref{refined-Strichartz} in Section 2 and Theorem \ref{multilinear-refined-strichartz} in Section 3. In Section 4, we prove Theorem \ref{Maximal-Schrodinger} using a bilinear argument. Section 5 and 6 are appendices concerning wave-packet decomposition and an epsilon removal lemma.

\medskip

{\bf Notations: }Throughout the paper, we will use the following notations:
\begin{enumerate}
    \item[$\bullet$] We use $a\sim b$ to mean that $ca\leq b\leq Ca$ for some unimportant constants $c$ and $C$ depending only on dimension. The symbol $B^n(x,r)$ represents the open ball centered at $x$, of radius $r$, in $\ZR^n$, and $B^n_r$ represents $B^n(0,r)$. $M,N,C$ are (big) constants depend only on dimension.
    \item[$\bullet$] For a rectangle $\om\in\ZR^n$, we use $c(\om)$ to denote the center of $\om$ and use $c(\om)_j$ to denote the $j$-th coordinate of $c(\om)$. 
    \item[$\bullet$] We write $A(R) \leq \mathrm{RapDec}(R)B$ to mean that for any power $N$, there is a constant $C_{N}$ such that
\begin{equation}
\nonumber
    A(R) \leq C_{N}R^{-N}B \;\; \text{for all $R \geq 1$}.
\end{equation}
    \item[$\bullet$] We set $w(x)=(1+|x|)^{-1000n}$ and $w_B(x)=w(\frac{x-c(B)}{R})$ for any $R$-ball $B\in\ZR^n$. We let $\|f\|_{L^p(w_B)}=\|fw_B\|_{L^p(\ZR^n)}$.
\end{enumerate}

\noindent
{\bf Acknowlegement}. I am deeply grateful to my advisor and teacher Prof. Xiaochun Li for introducing the problem to me and being patient and supportive all the time.

\section{Linear refined Strichartz estimates}
In this section we prove Theorem \ref{refined-Strichartz}. Similar to the argument in \cite{Du-Li-Guth-Schrodinger-R3} and \cite{Du-Li-Guth-Zhang-Schrodinger-Rn}, we need the following $l^2$ decoupling theorem for hypersurfaces with nonzero Gaussian curvature:
\begin{theorem}[Bourgain-Demeter]
\label{decoupling-theorem}
Consider the quadratic hypersurfaces in $\ZR^{n+1}$
\begin{equation}
\nonumber
    H^n_v=\{(\xi_1,\ldots,\xi_n,v_1\xi_1^2+\cdots+v_n\xi_n^2),|\xi|\leq1, v_j=\pm1\}.
\end{equation}
Let $\cn_\de(H^n_v)$ be a $\de$ neighbourhood of $H^n_v$ and let $\cp_\de$ be a finitely overlapping cover of $\cn_\de(H^n_v)$ with rectangles $\theta$ of dimensions $\de^{1/2}\times\cdots\times\de^{1/2}\times\de$. We denote by $f_\theta$ the Fourier restriction of $f$ to $\theta$, namely, $\wh{f_\theta}(\xi)=\wh{f}(\xi)\Id_\theta(\xi)$. Assuming ${\rm supp}(\wh{f})\subset\cn_\de$, then for any $\e>0$, $p\geq\frac{2(n+2-d(H_v^n))}{n-d(H_v^n)}$, we have
\begin{equation}
\label{decoupling}
    \|f\|_p\leq C_\e\de^{-\e}\de^{-\frac{n}{4}+\frac{n+2}{2p}}\Big(\sum_\theta\|f_\theta\|_p^2\Big)^{1/2}.
\end{equation}
\end{theorem}
The $l^2$ decoupling theorem was first proved by Bourgain and Demeter in \cite{Bourgain-Demeter-l2} for elliptic hypersurfaces. Then it was extended to general hypersurfaces with nonzero Gaussian curvature by the same authors in \cite{Bourgain-Demeter-hyperbolic-decoupling}. Via a simple induction argument (see \cite{Bourgain-Demeter-l2} Chapter 7), we have 
\begin{corollary}
Let $S\subset\ZR^{n+1}$ be a smooth compact hypersurface with nonzero Gaussian curvature. Let $\cn_\de(S),\cp_\de,f$ be as in Theorem \ref{decoupling-theorem}. Then when $p\geq\frac{2(n+2-d(S))}{n-d(S)}$
\begin{equation}
\label{decoupling1-2}
    \|f\|_p\leq C_\e\de^{-\e}\de^{-\frac{n}{4}+\frac{n+2}{2p}}\Big(\sum_\theta\|f_\theta\|_p^2\Big)^{1/2}.
\end{equation}
\end{corollary}
Here is a standard trick when we want to use decoupling for local $L^p$ norm $\|f\|_{L^p(B)}$, where $B$ is any $\de^{-1}$ ball in $\ZR^{n+1}$. Notice that we can replace $f_\theta$ by $(\wh{f}\wh\Id^\ast_\theta)^\vee$ for some bump function $\wh\Id^\ast_\theta$ on the right side of \eqref{decoupling1-2}. We choose a smooth cutoff function $\Id^\ast_B$ so that $\Id^\ast_B\gtrsim1$ on $B$ and $\wh\Id^\ast_B$ is supported in $B^{n+1}(0,\de)$. Then the Fourier transform of $f\Id^\ast_B$ is contained in $\cn_{2\de}(S)$. Hence we can apply \eqref{decoupling1-2} with $f_\theta$ replaced by $(\wh{f}\Id^\ast_\theta)^\vee$ to have
\begin{equation}
    \|f\Id_B^\ast\|_p\leq C_\e\de^{-\e}\de^{-\frac{n}{4}+\frac{n+2}{2p}}\Big(\sum_\theta\|\Id^\ast_\theta\ast(f\Id_B^\ast)\|_p^2\Big)^{1/2}.
\end{equation}
Since $\Id^\ast_\theta\ast(f\Id_B^\ast)=\Id^\ast_\theta\ast(f_{10\theta}\Id_B^\ast)$, we can use Hausdorff-Young inequality and the triangle inequality to bound
\begin{equation}
    \sum_\theta\|\Id^\ast_\theta\ast(f\Id_B^\ast)\|_p\lesssim\sum_\theta\|f_\theta\Id_B^\ast\|_p.
\end{equation}
Finally, since $\Id_B\lesssim \Id_B^\ast$ and since $\Id_B^\ast\lesssim w_B$, we obtain the local version of the decoupling estimate 
\begin{equation}
\label{decoupling2}
    \|f\|_{L^p(B)}\leq C_\e\de^{-\e}\de^{-\frac{n}{4}+\frac{n+2}{2p}}\Big(\sum_\theta\|f_\theta\|_{L^p(w_B)}^2\Big)^{1/2}.
\end{equation}
This is what we will use later in Section 2 and 3.

Unlike the situation in \cite{Du-Li-Guth-Zhang-Schrodinger-Rn}, we are facing a loss of $\de^{-\frac{n}{4}+\frac{n+2}{2p}}$ when applying the decoupling inequalities \eqref{decoupling2}. However, the loss can be picked up from quadratic (parabolic) rescaling. Henceforth, we believe the refined Strichartz estimates \eqref{refined-Strichartz-esti} remain true for $p>(2n+4)/n$, where $(2n+4)/n$ is the endpoint for Strichartz estimates for general surfaces of nonzero Gaussian curvature.

\begin{proof}[Proof of Theorem \ref{refined-Strichartz}]
We use induction on scale here. Our target scale is $R^{1/2}$. That is, we assume \eqref{refined-Strichartz-esti} holds for the scale $R^{1/2}$ and aim to prove \eqref{refined-Strichartz-esti} on the scale $R$.

We break the unit ball $B^n(0,1)$ in frequency space into finitely overlapping $R^{-1/4}$-cubes $\{q\}$, and break the physical space $\ZR^n$ into finitely overlapping  $R^{3/4}$-cubes $\{Q\}$. For both the frequency cover $\{q\}$ and the physical cover $\{Q\}$, we let $\{\Id_q^\ast\}$ and $\{\Id_Q^\ast\}$ be two smooth partitions of unity adaptd to them respectively. Now we define $f_{q,Q}=(\Id^\ast_q(f\Id_Q^\ast)^\wedge)^\vee$. Then heuristically, when restricting $\ce_{\Phi} f_{q,Q}$ to the ball $B^{n+1}(0,1)$, it is essentially supported in a fat tube $R_{q,Q}$ of dimensions $\sim R^{3/4}\times\cdots\times R^{3/4}\times R$, pointing to the direction $(\nabla\Phi(c(q)),-1)$. Let us put this heuristic claim into a lemma
\begin{lemma}
\label{support-lem-section-2}
Let $R_{q,Q}$ be a fat tube of radius $C R^{3/4}$ and length $CR$, pointing to direction $(\nabla\Phi(c(q)),-1)$, so that $10Q\subset R_{q,Q}\cap\{t=0\}$. Then for any $x\in B^{n+1}(0,R)$,
\begin{equation}
    \ce_\Phi f_{q,Q}(x)=(\Id_{R_{q,Q}}\ce_\Phi f_{q,Q})(x)+\rapid(R)\|f\|_2.
\end{equation}
\end{lemma}
\noindent Lemma \ref{support-lem-section-2} follows standardly via the method of (non) stationary phase. We omit details here.

In order to use the induction hypothesis \eqref{refined-Strichartz-esti} at the scale $R^{1/2}$ for $\ce_\Phi f_{q,Q}$, we need some preparations. First, we consider a collection of finitely overlapping rectangular tubes $S_j$ in $R_{q,Q}$ of dimensions $R^{1/2}\times\cdots\times R^{1/2}\times R^{3/4}$ and  having the same direction as $R_{q,Q}$. Notice that, after quadratic (parabolic) rescaling at the scale $R^{1/4}$, $S_j$ will become a $R^{1/4}$-cube and $R_{q,Q}$ will become a $R^{1/2}$-ball. We choose a smooth partition of unity $\{\Id^\ast_{S_j}\}$ associated to the collection of rectangular tube $\{S_j\}$, such that:
\begin{enumerate}
    \item $\Id_{S_j}^\ast\geq 0$ and $\Id^\ast_{S_j}(x)=\rapid(R)$ when $x\in \ZR^{n+1}\setminus 2S_j$.
    \item $\wh\Id^\ast_{S_j}\geq0$ and $\wh\Id^\ast_{S_j}$ is supported at a $R^{-1/2}\times\cdots\times R^{-1/2}\times R^{-3/4}$ rectangular tube $\om_{q,Q}$ centered at the origin, whose shortest side is parallel the longest side of $S_j$.
\end{enumerate}
To be able to use \eqref{refined-Strichartz-esti} as an induction hypothesis on the pair $(\{S_j\},R_{q,Q})$, we need a dyadic pigeonholing argument:
\begin{enumerate}
    \item Let $\la$ be a dyadic number, we sort the $S_j$ inside $2R_{q,Q}$ according to the quantity $\|(\ce_{\Phi} f_{q,Q})\Id_{2S_j}\|_p$. We define $\ZS_\la$ to be the collection of $S_j$ such that $\|(\ce_{\Phi} f_{q,Q})\Id_{2S_j}\|_p\sim\la$.
    \item For each $\la$, we sort the rectanglular tubes $S_j\in\ZS_\la$ by the number of such rectangles in  horizontal slabs of thickness $R^{3/4}$, which are perpendicular to the direction of $S_j$. Let $\mu$ be another dyadic number. We define $\ZS_{\la,\mu}$ to be the collection of $S_j\in\ZS_\la$ so that the number of $S_j$ in every horizontal slab is $\sim\mu$.
\end{enumerate}
We let $Y_{q,Q,\la,\mu}$ be the union of rectangles $2S_j$ with $S_j\in\ZS_{\la,\mu}$. For a fixed pair $(q,Q)$, we define $Y_{q,Q,\la,\mu}^\ast$ as
\begin{equation}
\label{Y-star}
    Y_{q,Q,\la,\mu}^\ast:=\sum_{S_j\in\ZS_{\la,\mu}}\Id_{S_j}^\ast
\end{equation}
Note that $\Id_{R_{q,Q}}(x)-\sum_{\la,\mu}Y_{q,Q,\la,\mu}^\ast(x)=\rapid(R)$ when $x\in R_{q,Q}$. Via Lemma \ref{support-lem-section-2}, for any $x\in B^{n+1}(0,R)$, one has
\begin{equation}
    \ce_{\Phi} f(x)=\ \sum_{\la,\mu}\Bigg(\sum_{q,Q}(\ce_{\Phi} f_{q,Q})Y_{q,Q,\la,\mu}^\ast\Bigg)(x)+\rapid(R)\|f\|_2.
\end{equation}
Since there are $\lesssim O(\log R)^2$ many dyadic numbers $\la,\mu$, by pigeonholing we can choose some particular $\la,\mu$ such that, there exists a set $Y'\subset Y$, $|Y|/|Y'|=O((\log R)^2)$ and for any $B_j\in Y'$,
\begin{equation}
\label{dyadic-preserving1}
    \|\ce_{\Phi} f\|_{L^p(B_j)}\lesssim (\log R)^2\Big\|\sum_{q,Q}(\ce_{\Phi} f_{q,Q})Y_{q,Q,\la,\mu}^\ast\Big\|_{L^p(B_j)}+\rapid(R)\|f\|_2.
\end{equation}

Next, we fix the dyadic numbers $\la,\mu$ in the rest of the argument, and abbreviate $Y_{q,Q,\la,\mu}(\text{resp.},Y_{q,Q,\la,\mu}^\ast)$ to $Y_{q,Q}(\text{resp.},Y_{q,Q}^\ast)$. Notice that the set $Y_{q,Q}$ we chose has a similar pattern to the $Y$ in \eqref{refined-Strichartz-esti}, with $\si\sim\si_{q,Q}$ being the value $\mu$ we have fixed. Since we only need to consider those $f_{q,Q}$ with $\|f_{q,Q}\|_2\gtrsim \rapid(R)\|f\|_2$, there are $O(\log R)$ possible dyadic values for these $\|f_{q,Q}\|_2$. Also, notice that for each $B_j\in Y'$ and a fixed frequency cube $q$, there is $O(1)$ physical cube $Q$ such that $B_j\subset Y_{q,Q}$. It implies that there are $O(\log R)$ possible dyadic values $\eta$ such that $\#\{(q,Q):B_j\subset Y_{q,Q}\}\sim\eta$. Therefore, by a dyadic pigeonholing argument, we can choose a dyadic value $\eta$, a set $\ZP_\eta$ of pairs $(q,Q)$ and a set $Y''\subset Y$ such that $|Y|/|Y''|=O((\log R)^4)$ and
\begin{enumerate}
    \item For all $(q,Q)\in\ZP_\eta$, $\|f_{q,Q}\|_2$ are the same up to a constant factor.
    \item For each $B_j\in Y''$, $\#\{(q,Q):B_j\subset Y_{q,Q}\}\sim\eta$.
    \item For each $B_j\in Y''$,
    \begin{equation}
\label{dyadic-preserving2}
    \|\ce_{\Phi} f\|_{L^p(B_j)}\lesssim (\log R)^4\Big\|\sum_{(q,Q)\in\ZP_\eta}(\ce_{\Phi} f_{q,Q})Y_{q,Q}^\ast\Big\|_{L^p(B_j)}+\rapid(R)\|f\|_2.
\end{equation}
\end{enumerate}
We fix $\eta$ in the rest of the argument.

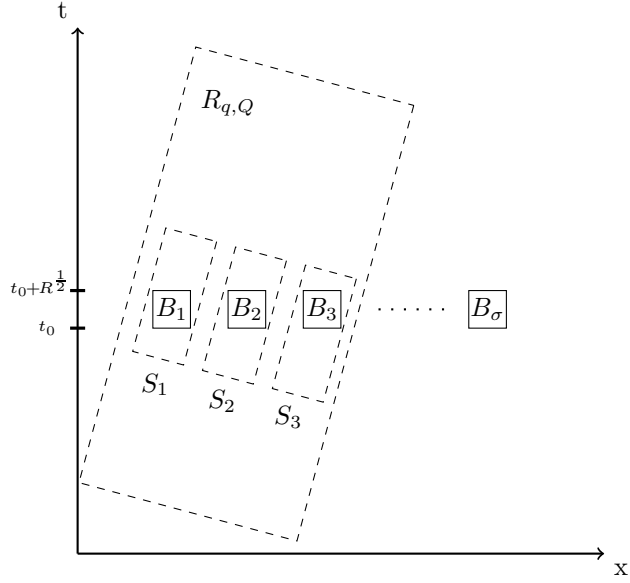
\begin{figure}
\begin{tikzpicture}
\draw[thick,->] (0,0) -- (7,0) node [anchor=north west] {x};
\draw[thick,->] (0,0) -- (0,7) node [anchor=south east] {t};
  \draw (1,3) -- ++(.5,0) -- ++(0,.5) -- ++(-.5,0) -- cycle;
  \draw (2,3) -- ++(.5,0) -- ++(0,.5) -- ++(-.5,0) -- cycle;
  \draw (3,3) -- ++(.5,0) -- ++(0,.5) -- ++(-.5,0) -- cycle;
  \draw (5.2,3) -- ++(.5,0) -- ++(0,.5) -- ++(-.5,0) -- cycle;
  \draw [loosely dotted, thick] (4,3.25) -- (5,3.25);
  \draw[dashed, rotate around={-15:(-3,9)}] (5,8) rectangle (2,2);
  \draw[dashed, rotate around={-15:(3,8.45)}]  (3,4) rectangle (2.3,2.3);
  \draw[dashed, rotate around={-15:(2.5,4.8)}] (3,4) rectangle (2.3,2.3);
  \draw[dashed, rotate around={-15:(2.05,1.15)}] (3,4) rectangle (2.3,2.3);
  \node at (2,6) {$R_{q,Q}$};
  \node at (1.25,3.25) {$B_1$};
  \node at (2.25,3.25) {$B_2$};
  \node at (3.25,3.25) {$B_3$};
  \node at (5.45,3.25) {$B_\si$};
  \node at (1.03,2.25) {$S_1$};
  \node at (1.92,2.03) {$S_2$};
  \node at (2.80,1.82) {$S_3$};
  \draw [very thick] (-0.1,3) -- (0.1,3);
  \draw [very thick] (-0.1,3.5) -- (0.1,3.5);
  \node at (-.4,3) {\tiny${t_0}$};
  \node at (-.5,3.6) {\tiny ${t_0\!\!+\!\!R^{\frac{1}{2}}}$};
\end{tikzpicture}
\caption{Possible position for $S_j$ and $B_j$.}
\label{figure}
\end{figure}

The following geometric fact plays a crucial role in the proof: Each $2S_j\in Y_{q,Q}$  contains $O(1)$ many cubes $B_j$ in a horizontal slab $\ZR^n\times\{t_0,t_0+R^{1/2}\}$. See Figure \ref{figure} for a possible position between $S_j$ and $B_j$. This is because the direction of $R_{q,Q}$ is $(\nabla\Phi(c(q)),-1)$, which makes an angle of $\sim1$ with respect to the horizontal plane $t=0$. Thus, 

\begin{equation}
    \frac{|Y_{q,Q}\cap Y|}{|Y|}\leq\frac{\si_{q,Q}}{\si}.
\end{equation}
Since $Y''\subset Y$, we have
\begin{equation}
\label{geometry-esti}
    \frac{|Y_{q,Q}\cap Y''|}{|Y''|}\leq(\log R)^4\frac{\si_{q,Q}}{\si}.
\end{equation}

Finally, combining the fact that for $B_j\in Y''$, $\#\{(q,Q):B_j\subset Y_{q,Q}\}\sim\eta$ and \eqref{geometry-esti}, we conclude

\begin{equation}
\label{geometry-esti2}
    \eta\leq (\log R)^4\frac{\si_{q,Q}}{\si}|\ZP_\eta|.
\end{equation}

\medskip

Now we are ready to begin our proof. Note that the Fourier transform of $(\ce_{\Phi} f_{q,Q})Y_{q,Q}^\ast$ is contained $N_{R^{-1/2}}(\Si)$ with $\Si=\{(\xi,\Phi(\xi)\}$. We use \eqref{decoupling2} for $\de=R^{-1/2}$ and \eqref{dyadic-preserving2} so that for each $B_j\in Y''$,
\begin{align}
    \|\ce_{\Phi} f\|_{L^p(B_j)} \lesssim&\, (\log R)^4\Bigg\|\sum_{(q,Q)\in\ZP_\eta}(\ce_{\Phi} f_{q,Q})Y_{q,Q}^\ast\Bigg\|_{L^p(B_j)}\!\!\!\!+\rapid(R)\|f\|_2\\ \nonumber
    \lesssim&\, R^{\frac{\e}{4}} R^{\frac{n}{8}-\frac{n+2}{4p}}\Bigg(\sum_{(q,Q)\in\ZP_\eta}\big\|(\ce_{\Phi} f_{q,Q})Y_{q,Q}^\ast\big\|_{L^p(w_{B_j})}^2\Bigg)^{1/2}\\[1ex] \nonumber
    &+\rapid(R)\|f\|_2.
\end{align}
Recall that the number of sets $Y_{q,Q}$ containing $B_j$ is $\sim\eta$. Since $Y_{q,Q}^\ast$ decays rapidly outside $Y_{q,Q}$, and since for fixed $q$, the set $Y_{q,Q}$ is finitely overlapped, we can apply H\"older's inequality to have
\begin{align}
\nonumber
    \Bigg\|\!\sum_{(q,Q)\in\ZP_\eta}(\ce_{\Phi} f_{q,Q})Y_{q,Q}^\ast\Bigg\|_{L^p(B_j)}\!\!\!\!\!\!\!\lesssim &\, \eta^{\frac{p-2}{2p}}R^{\frac{\e}{4}} R^{\frac{n}{8}-\frac{n+2}{4p}}\Bigg(\sum_{(q,Q)\in\ZP_\eta}\!\!\!\big\|(\ce_{\Phi} f_{q,Q})Y_{q,Q}^\ast\big\|_{L^p(w_{B_j})}^p\!\!\Bigg)^{1/p}\\ \nonumber
    &+\rapid(R)\|f\|_2.
\end{align}
Noticing that $\sum_{j}|w_{B_j}(x)|\lesssim1$, we raise $p$-th power to both sides and sum over all the $B_j\in Y''$ to get
\begin{equation}
\nonumber
    \big\|\ce_{\Phi} f\big\|_{L^p(Y'')}^p\lesssim \eta^{\frac{p-2}{2}}R^{\frac{\e p}{4}} R^{\frac{np}{8}-\frac{n+2}{4}}\!\!\!\!\sum_{(q,Q)\in\ZP_\eta}\!\!\!\big\|(\ce_{\Phi} f_{q,Q})Y_{q,Q}^\ast\big\|_p^p+\rapid(R)\|f\|_2^p.
\end{equation}
Since $|Y|\lesssim(\log R)^4|Y''|$, and since for all $B_j\in Y$, $\|\ce_{\Phi} f\|_{L^p(B_j)}$ are essentially the same, we have
\begin{equation}
\label{before-induction}
    \big\|\ce_{\Phi} f\big\|_{L^p(Y)}^p\lesssim \eta^{\frac{p-2}{2}}R^{\frac{\e p}{3}} R^{\frac{np}{8}-\frac{n+2}{4}}\!\!\!\!\sum_{(q,Q)\in\ZP_\eta}\!\!\!\big\|(\ce_{\Phi} f_{q,Q})Y_{q,Q}^\ast\big\|_p^p+\rapid(R)\|f\|_2^p.
\end{equation}

Let us take a close look at the term $\|(\ce_{\Phi} f_{q,Q})Y_{q,Q}^\ast\|_p^p$. Recall \eqref{Y-star} and recall the set $\ZS_{\la,\mu}$ from the paragraph above \eqref{Y-star}. Since the rectangular tubes $S_j$ are finitely overlapped and since $\Id_{S_j}^\ast$ decays rapidly outside $2S_j$, we have
\begin{align}
\label{smooth-to-sharp}
    \|(\ce_{\Phi} f_{q,Q})Y_{q,Q}^\ast\|_p^p&\lesssim\sum_{S_j\in\ZS_{\la,\mu}}\|(\ce_{\Phi} f_{q,Q})\Id_{S_j}^\ast\|_p^p\\ \nonumber
    &\lesssim \sum_{S_j\in\ZS_{\la,\mu}}\|(\ce_{\Phi} f_{q,Q})\Id_{2S_j}\|_p^p+\rapid(R)\|f\|_2^p.
\end{align}
Recall that $\|(\ce_{\Phi} f_{q,Q})\Id_{2S_j}\|_p$ are the same up to a constant, and recall that the rectangular tubes $S_j\in\ZS_{\la,\mu}$ are arranged in horizontal slabs of thickness $R^{3/4}$ which are perpendicular to the direction of $S_j$, so that each of these horizontal slabs contains $\sim\mu=\si_{q,Q}$ many $S_j$ (See Figure \ref{figure} again). Hence the doubles $2S_j$ are arranged in the same way, and each thicker horizontal slab contains $\sim\si_{q,Q}$ many doubles $2S_j$. Note that the doubles $2S_j$ are finitely overlapped. Thus, at a cost of an absolute constant, we can use the induction hypothesis \eqref{refined-Strichartz-esti} at scale $R^{1/2}$ after quadratic (parabolic) rescaling (See \eqref{scaling}, \eqref{scaling2} and \eqref{scaling3}) to obtain
\begin{equation}
\label{scale-half}
    \sum_{S_j\in\ZS_{\la,\mu}}\|(\ce_{\Phi} f_{q,Q})\Id_{2S_j}\|_p^p\lesssim R^{\e p/2}\si_{q,Q}^{-\frac{p-2}{2}}\|f_{q,Q}\|_2^p
\end{equation}

Combine \eqref{before-induction}, \eqref{smooth-to-sharp} and \eqref{scale-half} so that
\begin{equation}
\label{local-estimate}
    \big\|\ce_{\Phi} f\big\|_{L^p(Y)}^p\lesssim\eta^{\frac{p-2}{2}}\si_{q,Q}^{-\frac{p-2}{2}}R^{\frac{3\e p}{4}}\sum_{(q,Q)\in\ZP_\eta}\big\|f_{q,Q}\big\|_2^p+\rapid(R)\|f\|_2^p.
\end{equation}
Now we can apply the geometric estimate \eqref{geometry-esti2} to have
\begin{equation}
    \big\|\ce_{\Phi} f\big\|_{L^p(Y)}^p\lesssim\si^{-\frac{p-2}{2}}R^{\e p}|\ZP_\eta|^{\frac{p-2}{2}}\sum_{(q,Q)\in\ZP_\eta}\big\|f_{q,Q}\big\|_2^p+\rapid(R)\|f\|_2^p.
\end{equation}
Since $\|f_{q,Q}\|_2$ are the same up to a constant, we get
\begin{align}
    \big\|\ce_{\Phi} f\big\|_{L^p(Y)}^p&\lesssim\si^{-\frac{p-2}{2}}R^{\e p}\Bigg(\sum_{(q,Q)\in\ZP_\eta}\big\|f_{q,Q}\big\|_2^2\Bigg)^{p/2}+\rapid(R)\|f\|_2^p\\ \nonumber
    &\lesssim\si^{-\frac{p-2}{2}}R^{\e p}\|f\|_2^p.
\end{align}
Finally, we take $p$-th root to both sides to have \eqref{refined-Strichartz-esti} for the scale $R$, and hence close the induction. 
\end{proof}

\begin{remark}
\rm

Although it is shown in \cite{Du-Li-Guth-Zhang-Schrodinger-Rn} that the refined Strichartz estimates \eqref{refined-Strichartz-esti} are sharp with respect to $\si$ when $p=(2n+4)/n$, we do not know whether it is sharp for $p>(2n+4)/n$. The example in \cite{Du-Li-Guth-Zhang-Schrodinger-Rn} does not work as a counterexample here. For $p>(2n+4)/n$, we may gain more from the "spreading out" property \eqref{spreading-out-condition} (or maybe a refined one).
\end{remark}

\section{Multilinear analogue}

This section is devoted to prove the multilinear refined Strichartz estimates \eqref{multilienar-refined-strichartz-esti}. Similar to the argument in \cite{Du-Li-Guth-Zhang-Schrodinger-Rn}, we will need the following multilinear Kakeya estimates from \cite{BCT-multilinear}. See also \cite{Guth-endpoint} and \cite{Guth-multiliear}.
\begin{theorem}
Suppose $S_j\subset S^{n-1},j=1,2,\ldots,k$, $k\geq2$. Suppose $l_{j,a}$ are lines in $\ZR^n$ and the direction of $l_{j,a}$ lies in $S_j$. Suppose that for any $v_j\in S_j$, we have
\begin{equation}
    |v_1\wedge\cdots\wedge v_k|\geq c>0.
\end{equation}
Let $T_{j,a}$ be the characteristic function for the $1$-neighbourhood of $l_{j,a}$. Then for any $\e>0$, we have
\begin{equation}
\label{multilinear-Kakeya-esti}
    \int_{B^n_R}\prod_{j=1}^k\Big(\sum_{a=1}^{N_j}T_{j,a}\Big)^{1/(k-1)}\leq C_\e R^\e\prod_{j=1}^kN_j^{1/(k-1)}.
\end{equation}
\end{theorem}
\begin{proof}[Proof for Theorem \ref{multilinear-refined-strichartz}.]
For all $j$, we process $\|\ce_{\Phi} f_j\|_{\lp(Y)}$ following the proof in Theorem \ref{refined-Strichartz} simultaneously to have $\la_j,\mu_j,\eta_j$, a set $\ZP_{\eta_j}$ and a set $Y''\subset Y$. We abbreviate $\ZP_{\eta_j}$ to $\ZP_j$ in the rest of the argument. Thus, similar to \eqref{local-estimate} one has
\begin{equation}
\label{priori-estimate}
    \big\|\ce_{\Phi} f_j\big\|_{L^p(Y)}\lesssim\eta_j^{\frac{p-2}{2p}}\si_{j,q,Q}^{-\frac{p-2}{2p}}R^{\frac{3\e }{4}}\sum_{(q,Q)\in\ZP_j}\big\|f_{q,Q}\big\|_2+\rapid(R)\|f_j\|_2.
\end{equation}
The set $Y''$ was chosen such that $|Y|\lesssim(\log R)^{4k}|Y''|$, and for each $R^{1/2}$ cube $B\in Y''$ and each $j$,
\begin{equation}
    \#\{(q,Q)\in \ZP_j:B\subset Y_{q,Q}\}\sim\eta_j.
\end{equation}
Note that the set $Y$ is a union of $N$ many $R^{1/2}$ cubes. Therefore, summing all the $B\subset Y''$ in the above estimate we have
\begin{equation}
\label{simple-fact}
    N\prod_{j=1}^k\eta_j^{\frac{1}{k-1}}\leq(\log R)^{4k}\sum_{B\in Y''}\prod_{j=1}^k\big(\#\{(q,Q)\in \ZP_j:B\subset Y_{q,Q}\}\big)^{\frac{1}{k-1}}.
\end{equation}

Our goal is to show
\begin{equation}
\label{kakeya-type-estimate}
    N\prod_{j=1}^k\eta_j^{\frac{1}{k-1}}\leq C_\e R^\e\prod_{j=1}^k\big(\si_{j,q,Q}|\ZP_j|\big)^{\frac{1}{k-1}}.
\end{equation}
Since $Y_{q,Q}$ is a collection of tubes that each tube has dimensions $R^{1/2}\times\cdots\times R^{1/2}\times R^{3/4}$, we are motivated to use $k$-linear Kakeya estimate \eqref{multilinear-Kakeya-esti} on the scale $R^{3/4}$. We break the ball $B^{n+1}_R$ into lattice $R^{3/4}$-cubes $K$. Define
\begin{equation}
    \ZP_{j,K}=\{(q,Q)\in\ZP_j:K\subset 2R_{q,Q}\}.
\end{equation}
Then \eqref{kakeya-type-estimate} is equivalent to
\begin{equation}
\label{Kakeya-type-estimate2}
    N\prod_{j=1}^k\eta_j^{\frac{1}{k-1}}\lesssim R^{\e}\sum_K\sum_{B\in Y'',B\subset K}\prod_{j=1}^k\big(\#\{(q,Q)\in \ZP_{j,K}:B\subset Y_{q,Q}\}\big)^{\frac{1}{k-1}}.
\end{equation}
On each in $K$, we have $k$ transverse collections of tubes of dimensions $R^{1/2}\times\cdots\times R^{1/2}\times R^{3/4}$ passing through it, and the number of such tubes in the $j$-th collection is $\lesssim|\ZP_{j,K}|\cdot\si_{j,q,Q}$. Since the intersection of the tubes in different collections is essentially a $R^{1/2}$ cube, we can use $k$-linear Kakeya estimate \eqref{multilinear-Kakeya-esti} to have
\begin{equation}
\nonumber
    \sum_{B\in Y'',B\subset K}\prod_{j=1}^k\big(\#\{(q,Q)\in \ZP_{j,K}:B\subset Y_{q,Q}\}\big)^{\frac{1}{k-1}}\lesssim R^\e\prod_{i=1}^k\big(|\ZP_{j,K}|\cdot\si_{j,q,Q}\big)^\frac{1}{k-1}.
\end{equation}
Plug this back to \eqref{Kakeya-type-estimate2} so that
\begin{equation}
    N\prod_{j=1}^k\eta_j^{\frac{1}{k-1}}\lesssim R^\e\sum_{K} \prod_{i=1}^k\big(|\ZP_{j,K}|\cdot\si_{j,q,Q}\big)^\frac{1}{k-1}.
\end{equation}
Notice that any $R_{q,Q}\in\ZP_{j,K}$ is in fact a tube of dimensions $R^{3/4}\times\cdots\times R^{3/4}\times R$. We thus can use $k$-linear Kakeya again to have
\begin{equation}
    \sum_{K} \prod_{i=1}^k\big(|\ZP_{j,K}|\cdot\si_{j,q,Q}\big)^\frac{1}{k-1}\lesssim R^{\e}\prod_{j=1}^k\big(|\ZP_j|\cdot\si_{j,q,Q}\big)^{\frac{1}{k-1}}.
\end{equation}
Hence we get \eqref{kakeya-type-estimate}.

To finish the proof, we first invoke H\"older's inequality so that
\begin{align}
    \Big\|\prod_{j=1}^k|\ce_{\Phi}f_j|^{1/k}\Big\|_{\lp(Y)}^k &\leq \prod_{j=1}^k\Big\|\ce_{\Phi}f_j\Big\|_{\lp(Y)}. 
\end{align}
Via \eqref{priori-estimate} and \eqref{kakeya-type-estimate}, we can bound the right hand side of the above estimate as
\begin{align}
\nonumber
    \prod_{j=1}^k\Big\|\ce_{\Phi}f_j\Big\|_{\lp(Y)}\!&\lesssim\Big(\prod_{j=1}^k\eta_j^{\frac{p-2}{2p}}\si_{j,q,Q}^{-\frac{p-2}{2p}}R^{\frac{3\e}{4}}\!\!\!\!\sum_{(q,Q)\in\ZP_j}\big\|f_{q,Q}\big\|_2+\rapid(R)\|f_j\|_2\Big)\\ \nonumber
    &\lesssim N^{-\frac{(k-1)(p-2)}{2p}}R^{\e}\Big(\prod_{j=1}^k|\ZP_j|^{\frac{p-2}{2p}}\!\!\!\!\sum_{(q,Q)\in\ZP_j}\big\|f_{q,Q}\big\|_2+\rapid(R)\|f_j\|_2\Big).
\end{align}
Note that all $(q,Q)\in\ZP_j$, $\|f_{q,Q}\|_2$ are the same up to a constant. Via the $L^2$ orthgonality among $f_{q,Q}$, the above two estimates yield 
\begin{equation}
    \Big\|\prod_{j=1}^k|\ce_{\Phi}f_j|^{1/k}\Big\|_{\lp(Y)}^k\lesssim N^{-\frac{(k-1)(p-2)}{2p}}R^{\e}\prod_{j=1}^k\|f_j\|_2.
\end{equation}
We take the $k$-th root on both sides to finish the proof of \eqref{multilienar-refined-strichartz-esti}.
\end{proof}

\section{Maximal extension operator}

We will prove Theorem \ref{Maximal-Schrodinger} here. The methods we use are the sharp $L^2$ estimate in \cite{Du-Zhang} and the bilinear argument in \cite{Tao-Bilinear}. Following an epsilon removal argument that we will prove in the Section 6, it suffices to prove a local result:
\begin{lemma}
\label{Maximal-Schrodinger-local}
Suppose that $\Phi:B^n(0,2)\to\ZR$ is a function satisfying $d(\Phi)=0$, $|\nabla\Phi|\lesssim1$ and $|\det(\nabla^2\Phi)|\sim1$. Let $f:\ZR^n\to\ZC$ be with $\wh{f}\in L^2(B^n(0,1))$. Then for any $p>2+\frac{4}{n+2-1/n}$ and any $\e>0$, there exists a constant $C_\e$ such that for all $R\geq 1$,
\begin{equation}
\label{Maximal-Schrodinger-local-esti}
    \big\|\sup_{|t|<R}|\ce_{\Phi}f|\big\|_{L^p(B^n_R)}\leq C_\e R^\e\|f\|_2.
\end{equation}
\end{lemma}


Since $\wh{f}$ is supported in the unit ball $B^n(0,1)$, we can freely add a bump function $\wh\phi$ that equals to $1$ in $B^n(0,1)$ to the definition of the extension operator $\ce_\Phi$ as
\begin{equation}
\label{extension-2}
    \ce_{\Phi} f(x,t)=\int_{\ZR^n}e^{i(x\cdot\xi+t\Phi(\xi))}\wh\phi(\xi)\wh{f}(\xi)d\xi.
\end{equation}
Here we also assume that $\wh\phi$ is supported in $B^n(0,2)$. Let us take \eqref{extension-2} as the definition of the extension operator $\ce_\Phi$ in the rest of this section. 

\medskip

To prove Lemma \ref{Maximal-Schrodinger-local}, we will do induction on the scale $R$ for smooth functions like $\wh\phi$ in \eqref{extension-2}, and for phase functions falling into the large family 
\begin{equation}
\label{family-of-functions}
    \Xi:=\{\Phi:B^n(0,2)\to\ZR,~|\nabla\Phi|\lesssim1~\text{and}~|\det(\nabla^2\Phi)|\sim1\}.
\end{equation}
Also, we need a bilinear version of Lemma \ref{Maximal-Schrodinger-local}, namely, 
\begin{lemma}
\label{bilinear-Maximal-Schrodinger}
Suppose that $\Phi:B^n(0,2)\to\ZR$ is a function satisfying $d(\Phi)=0$, $|\nabla\Phi|\lesssim1$ and $|\det(\nabla^2\Phi)|\sim1$ and suppose that $f_1,f_2$ have frequencies 2-transversely supported in $B^n(0,1)$. Then for $p>1+\frac{2}{n+2-1/n}$ and any $\e>0$, there exist a constant $C_\e$ that for all $R\geq 1$,
\begin{equation}
\label{bilinear-Maximal-Schrodinger-esti}
    \big\|\sup_{|t|<R}|\ce_{\Phi}f_1\ce_{\Phi}f_2|\big\|_{L^p(B^n_R)}\leq C_\e R^\e\|f_1\|_2\|f_2\|_2.
\end{equation}
\end{lemma}

First we will show how \eqref{bilinear-Maximal-Schrodinger-esti} implies \eqref{Maximal-Schrodinger-local-esti}. We follow the idea in \cite{Bourgain-Guth-Oscillatory} to break $\ce_\Phi f$ into linear and bilinear parts. For a large number $K>1$, let $\{\tau\}$ be a collection of finitely overlapping $K^{-1}$-cubes in $B^n(0,1)$, and let $\wh\vp_\tau(\xi):=\wh\vp((\xi-c(\tau))/K)$ be an associated smooth partition of unity to the unit ball $B^n(0,1)$ that $\wh\vp$ satisfies the following two properties:
\begin{enumerate}
    \item $\wh\vp$ is supported in $B^n(0,2)$.
    \item $|\partial^\al\wh\vp|\leq C_\al$ for every multi-index $\al\in\ZN^n$.
\end{enumerate}
We set $K=R^{\e^2}$ and define $f_\tau:=\vp_\tau\ast f$, so $f=\sum_\tau f_\tau$. Since for each $\tau$, there are only $\sim n$ many lattice $K^{-1}$-cubes that are adjacent to it, we have
\begin{equation}
    |\ce_{\Phi}f|\leq C\max_{\tau}|\ce_{\Phi}f_\tau|+\sum_{(\tau_1,\tau_2)\in\La_K}|\ce_{\Phi}f_\tau|^{\frac{1}{2}}\cdot|\ce_{\Phi}f_\tau|^{\frac{1}{2}},
\end{equation}
where $\La_K:=\{(\tau_1,\tau_2):{\rm dist}(\tau_1,\tau_2)\geq 2K^{-1}\}$. Take $L^p$ norm to both sides and use triangle inequality so that
\begin{align}
\label{dichotomy}
    \big\|\sup_{|t|<R}\ce_{\Phi}f|\big\|_{\lp(B^n_R)} \leq&\, C\big\|\sup_{|t|<R}\max_{\tau}|\ce_{\Phi}f_\tau|\big\|_{\lp(B^n_R)}\\ \label{bilinear-part}
    &+\sum_{(\tau_1,\tau_2)\in\La_K}\big\|\sup_{|t|<R}|\ce_{\Phi}f_{\tau_1}|^{\frac{1}{2}}|\ce_{\Phi}f_{\tau_2}|^{\frac{1}{2}}|\big\|_{\lp(B^n_R)}
\end{align}

For the linear part \eqref{dichotomy}, via the embedding $l^\infty\subset l^p$, one has
\begin{equation}
\label{lp-l-infty}
    \big\|\sup_{|t|<R}\max_{\tau}|\ce_{\Phi}f_\tau|\big\|_{\lp(B^n_R)}\leq\Bigg(\sum_\tau\big\|\sup_{|t|<R}|\ce_{\Phi}f_\tau|\big\|_{\lp(B^n_R)}^p\Bigg)^{1/p}.
\end{equation}
Let us fix a $K^{-1}$ cap $\tau$. Since $\wh{f_\tau}$ is supported in $2\tau$, we can freely add a cutoff function $\Id_{2\tau}$ to $\wh{f_\tau}$ and write
\begin{align}
\nonumber
    \big\|\sup_{|t|<R}|\ce_{\Phi}f_\tau|\big\|_{\lp(B^n_R)}^p&=\int_{B_R}\sup_{|t|\leq R}\Big|\int e^{ix\cdot\xi+it\Phi(\xi)}\wh\phi(\xi)\wh\vp\Big(\frac{\xi-c(\tau)}{K^{-1}}\Big)\wh{f}(\xi)d\xi\Big|^pdx\\
    &=\int_{B_R}\sup_{|t|\leq R}\Big|\int e^{ix\cdot\xi+it\Phi(\xi)}\wh\phi(\xi)\wh\vp\Big(\frac{\xi-c(\tau)}{K^{-1}}\Big)(\Id_{2\tau}\wh{f})(\xi)d\xi\Big|^pdx.
\end{align}
Consider the change of variable
\begin{equation}
\label{scaling}
    \eta:=\frac{\xi-{c(\tau)}}{K^{-1}}\,\,\, \Longleftrightarrow\,\,\, \xi=K^{-1}\eta+c(\tau)
\end{equation}
and a new function function $\Psi(\eta)$ defined as
\begin{equation}
\label{scaling2}
    \Psi(\eta):=K^{2}[\Phi({c(\tau)}+K^{-1}\eta)-\Phi({c(\tau)})-K^{-1}\nabla\Phi({c(\tau)})\cdot\eta].
\end{equation}
One can check directly that $\Psi$ falls into the family $\Xi$ introduced in \eqref{family-of-functions}. Now let us introduce the parabolic rescaling 
\begin{equation}
\label{scaling3}
    \cl_\tau:(y,u)=\big(K^{-1}(\nabla\Phi(c(\tau))+x),K^{-2}t\big), 
\end{equation}
so that for another smooth function $\wh\psi$ supported on $B^n(0,2)$, one has
\begin{align}
    &\int_{B_R}\sup_{|t|\leq R}\Big|\int e^{ix\cdot\xi+it\Phi(\xi)}\wh\phi(\xi)\wh\vp\Big(\frac{\xi-c(\tau)}{K^{-1}}\Big)(\Id_{2\tau}\wh{f})(\xi)d\xi\Big|^pdx\\ \nonumber
    =&\,K^{-np+n}\int_{B_{RK^{-1}}}\sup_{|u|\leq RK^{-2}}\Big|\int e^{iy\cdot\eta+iu\Psi(\xi)}\wh\psi(\eta)(\Id_{2\tau}\wh{f})(K^{-1}\eta+c(\tau))d\eta\Big|^pdy,
\end{align}
which, after using \eqref{Maximal-Schrodinger-local-esti} as an induction hypothesis at scale $RK^{-1}$, can be bounded above by
\begin{equation}
    C_\e^pK^{-\frac{pn}{2}+n}R^{p\e} K^{-p\e}\|\Id_{2\tau}\wh{f}\|_2^p.
\end{equation}
Hence we can sum up all the $K^{-1}$ caps $\tau$ and use Minkowski's inequality to obtain
\begin{align}
\label{part1}
    \Bigg(\sum_\tau\big\|\sup_{|t|<R}|\ce_{\Phi}f_\tau|\big\|_{\lp(B^n_R)}^p\Bigg)^{1/p}&\lesssim C_\e K^{-\frac{n}{2}+\frac{n}{p}-\e}R^\e\Big(\sum_\tau\|\Id_{2\tau}\wh{f}\|_2^2\Big)^{1/2}\\ \nonumber
    &\leq C_\e K^{-\frac{n}{2}+\frac{n}{p}-\e}R^\e\|f\|_2.
\end{align}

\medskip

Next, we will use the bilinear estimate \eqref{bilinear-Maximal-Schrodinger-esti} to bound the bilinear part \eqref{bilinear-part}. Invoking parabolic rescaling and an affine transformation when needed, we use \eqref{bilinear-Maximal-Schrodinger-esti} with $\e$ replaced by a smaller factor $\e^2$ to have
\begin{equation}
    \big\|\sup_{|t|<R}|\ce_{\Phi}f_{\tau_1}|^{\frac{1}{2}}|\ce_{\Phi}f_{\tau_2}|^{\frac{1}{2}}\big\|_{L^p(B^n_R)}\lesssim C_{\e^2}K^CR^{\e^2}\big\|f_{\tau_1}\big\|_2^{1/2}\big\|f_{\tau_2}\big\|_2^{1/2}.
\end{equation}
Sum up all the $(\tau_1,\tau_2)\in \La$ so that
\begin{equation}
\nonumber
    \sum_{(\tau_1,\tau_2)\in\La_K}\!\!\!\!\big\|\sup_{|t|<R}|\ce_{\Phi}f_{\tau_1}|^{\frac{1}{2}}|\ce_{\Phi}f_{\tau_2}|^{\frac{1}{2}}|\big\|_{\lp(B^n_R)}\leq C_{\e^2}K^CR^{\e^2}\!\!\!\!\sum_{(\tau_1,\tau_2)\in\La_K}\!\!\!\!\big\|f_{\tau_1}\big\|_2^{1/2}\big\|f_{\tau_2}\big\|_2^{1/2}.
\end{equation}
We apply H\"older's inequality to the right hand side of the above inequality to get
\begin{eqnarray}
\label{part2}
    \sum_{(\tau_1,\tau_2)\in\La_K}\!\!\!\!\!\!\big\|\sup_{|t|<R}|\ce_{\Phi}f_{\tau_1}|^{\frac{1}{2}}|\ce_{\Phi}f_{\tau_2}|^{\frac{1}{2}}|\big\|_{\lp(B^n_R)} \!\!\!\!\!\!&\leq&\!\!\!\!\! C_{\e^2}K^{2C}R^{\e^2}\big(\!\sum_{\tau_1,\tau_2}\!\|f_{\tau_1}\|_2^2\|f_{\tau_2}\|_2^2\big)^{1/4}\\ \nonumber
    &\leq&\!\!\!\!\! C_{\e^2}K^{2C}R^{\e^2}\|f\|_2.
\end{eqnarray}

Finally, we combine \eqref{dichotomy}, \eqref{bilinear-part}, \eqref{part1} and \eqref{part2} to have
\begin{equation}
     \big\|\sup_{|t|<R}|\ce_{\Phi}f|\big\|_{\lp(B^n_R)} \leq C_\e\big( CK^{-\frac{n}{2}+\frac{n}{p}-\e}R^\e+C_{\e^2}K^{2C}R^{\e^2}\big)\|f\|_2.
\end{equation}
The induction closes as 
\begin{equation}
\nonumber
    CK^{-\frac{n}{2}+\frac{n}{p}-\e}R^\e+C_{\e^2}K^{2C}R^{\e^2}=C R^{\e^2(-\frac{n}{2}+\frac{n}{p}-\e)}R^\e+C_{\e^2}R^{2C\e^2+\e^2}\leq R^{\e}. 
\end{equation}
when the radius $R$ is large enough and $p>2$. \qed

\medskip

It remains to show \eqref{bilinear-Maximal-Schrodinger-esti}. For convenience we make a definition about certain sets in $\ZR^{n+1}$:

\begin{definition}[Horizontally sparse]
\label{set-X}
Suppose that $X$ is a union of unit cubes in $\ZR^{n+1}$. We say $X$ is ``horizontally sparse", if for any lattice $1$-cube $U\subset B^n(0,R)$, there are only $O(1)$ many $\ZR^{n+1}$ unit cubes in $X$ whose intersection with $ P_U$ is not empty. Here the vertical stripe $P_U$ is defined as:
\begin{equation}
    P_U:=\{(x,t):x\in U, t\in\ZR\}.
\end{equation}
\end{definition}
Since $\wh{f}\subset B^n(0,1)$, $\ce_{\Phi}f$ is essentially constant on every 1-cube in $\ZR^{n+1}$. In fact, similar to Lemma \ref{weight}, we can find a horizontally sparse set $X\subset B^{n+1}(0,R)$ and a small factor $\be=\e^{1000}$ so that
\begin{align}
    \big\|\sup_{|t|<R}|\ce_{\Phi}f_1|^{\frac{1}{2}}|\ce_{\Phi}f_2|^{\frac{1}{2}}|\big\|_{\lp(B^n_R)} \lesssim&\, R^{O(\be)} \big\||\ce_{\Phi}f_1|^{\frac{1}{2}}|\ce_{\Phi}f_2|^{\frac{1}{2}}|\big\|_{\lp(X)}\\ \nonumber
    &+ \rapid(R)\|f_1\|_2^{1/2}\cdot\|f_2\|_2^{1/2}.
\end{align}
Let us fix the set $X$ from now on. Conversely, for an arbitrary $r$ ball $B\subset B^{n+1}(0,R)$ with $r\geq1$, one has
\begin{align}
\label{X-other-direction}
    \big\||\ce_{\Phi}f_1|^{\frac{1}{2}}|\ce_{\Phi}f_2|^{\frac{1}{2}}|\big\|_{\lp(X\cap B)}\lesssim \big\|\sup_{|t|<R}|\ce_{\Phi}f_1|^{\frac{1}{2}}|\ce_{\Phi}f_2|^{\frac{1}{2}}|\big\|_{\lp(B\cap\ZR^n)}.
\end{align}

\medskip

We will adapt the argument in \cite{Tao-Bilinear}. First we state a wave-packet decomposition and will prove it in the appendix. A similar formulation can be founded in \cite{Tao-Bilinear} and \cite{Guth-restriction-R3}. See also \cite{Lacey-Thiele-BHilbert}.

Let ${\bf q}=\{q\}$ be the collection of lattice $R^{-1/2}$-cubes in $B^n(0,1)$. Define $\ZT_q$ be a collection of rectangles $T$ in $B^{n+1}(0,2R)$ of dimensions $R^{1/2}\times\cdots\times R^{1/2}\times R$ such that the projection of $T$ onto the hyperplane orthogonal to the vector $(\nabla\Phi(c(q)),-1)$ is a lattice $R^{1/2}$-cube in $B^n(0,R)$. We let $\ZT=\bigcup_q\ZT_q$. Now we are ready to state our wave-packet decomposition.

\begin{proposition}
\label{wave-packet-paraboloid}
Assume $f\in L^2(B^n(0,1))$. Let $\de$ be a small number taking care of Schwartz tails. Then for each $f_q$, we can pick a collection of functions $f_T$, $T\in\ZT$ such that
\begin{enumerate}
    \item {\rm supp}$(\wh{f_T})\subset 2q$~{\rm for~some~}$q\in{\bf q}$.
    \item $|\ce_{\Phi} f_T(x,t)|\lesssim \rapid(R)\|f\|_2~for~(x,t)\in B^{n+1}_R\setminus R^{\de}T$.
    \item $\big|\ce_{\Phi} f(x,t)-\sum_{T\in\ZT}\ce_{\Phi} f_T(x,t)\big|\lesssim \rapid(R)\|f\|_2$~for $(x,t)\in B^{n+1}(0,R)$.
    \item $\sum_{T\in\ZT}\|f_T\|_2^2\sim\|f\|_2^2$.
    \item $\langle \ce_{\Phi} f_T,\ce_{\Phi} f_{T'}\rangle=\langle f_T ,f_{T'}\rangle=0, $~{\rm if}~{\rm~dist(}{\rm supp}$(\wh{f_T})$,~{\rm supp}$(\wh{f_{T'}}))\geq3R^{-1/2}$;\\
    $|\langle \ce_{\Phi} f_T,\ce_{\Phi} f_{T'}\rangle|\lesssim R|\langle f_T ,f_{T'}\rangle|\lesssim \rapid(R) $, {\rm if} {\rm dist($T,T'$)$\geq R^{1/2+\de}$.}
    \item $\|\ce_{\Phi}f_T\|_{L^2(B^{n+1}_R)}\lesssim R^{1/2}\|f_T\|_2$.
\end{enumerate}
We call $\sum_{T\in\ZT}f_T$ a wave packet decomposition of $f$.
\end{proposition}

Let $\sum_{T_j\in\ZT_j}f_{T_j}$ be a wave-packet decomposition of $f_j$, $j=1,2$. We set $\phi_T:=R^{-n/4}\chi_{B_{2R}^{n+1}}\ce_{\Phi} f_T$ be an $L^2$ normalized wave packet, so that
\begin{align}
    \big\|\sup_{|t|<R}|\ce_{\Phi}f_1|^{\frac{1}{2}}|\ce_{\Phi}f_2|^{\frac{1}{2}}|\big\|_{\lp(B^n_R)} \lesssim&\, R^{n/4}\big\|\sum_{T_1\in\ZT_1}\sum_{T_2\in\ZT_2}\phi_{T_1}\phi_{T_2}\big\|_{L^{p/2}(X)}^{1/2}\\ \nonumber
    &+ \rapid(R)\|f_1\|_2^{1/2}\cdot\|f_2\|_2^{1/2}.
\end{align}
Thus, it suffices to show that
\begin{equation}
    R^{n/4}\big\|\sum_{T_1\in\ZT_1}\sum_{T_2\in\ZT_2}\phi_{T_1}\phi_{T_2}\big\|_{L^{p}(X)}\leq C_\e R^\e\|f_1\|_2\|f\|_2.
\end{equation}
for any $\e>0$ and any $p>1+\frac{2}{n+2-1/n}$. 

Without loss of generality, we assume $\|f_1\|_2=\|f\|_2=1$, so $\|\phi_{T_j}\|_\infty\leq1$ . Since we only need to consider those $\phi_{T_j}$ with $\|\phi_{T_j}\|_\infty$ ranging in $[R^{-C},1]$, and there $\lesssim O(\log R)$ many dyadic value in $[R^{-C},1]$. Thus, we can assume that for each $j=1,2$, $\|\phi_{T_j}\|_\infty$ are equal up to a factor smaller than $2$. Invoking Property $(4)$ and $(5)$ in Proposition \ref{wave-packet-paraboloid}, it suffices to show
\begin{equation}
\label{kakeya-type}
    \big\|\sum_{T_1\in\ZT_1}\sum_{T_2\in\ZT_2}\phi_{T_1}\phi_{T_2}\big\|_{L^{p}(X)}\leq C_\e R^\e|\ZT_1|^{1/2}|\ZT_2|^{1/2}.
\end{equation}

Let ${\bf Q}=\{Q\}$ be the lattice $R^{1/2}$ cubes in $B^{n+1}(0,2R)$. We define $Q(T_j)=\{Q\in{\bf Q}:R^\de Q\cap T_j\not=\varnothing\}$ and $\ZT_j(Q)=\{T_j\in\ZT_j:R^\de Q\cap T_j\not=\varnothing\}$. Since $1\leq |{\bf Q}|,|\ZT_j(Q)|,|Q(T_j)|\leq R^N$, by the same dyadic argument above, we only need to show that for every pair of dyadic values $(\la_j,\mu_j),\la_j,\mu_j\in[1,R^N]$,
\begin{equation}
    \big\|\sum_{Q\in{\bf Q}[\mu_1,\mu_2]}\sum_{T_j\in\ZT_j[\la_j,\mu_1,\mu_2]}\phi_{T_1}\phi_{T_2}\Id_Q\big\|_{L^{p}(X)}\leq C_\e R^\e|\ZT_1|^{1/2}|\ZT_2|^{1/2}.
\end{equation}
where ${\bf Q}[\mu_1,\mu_2]=\{Q\in{\bf Q},~|\ZT_j(Q)|\sim\mu_j,j=1,2\}$ and $\ZT_j[\la_j,\mu_1,\mu_2]=\{T_j\in\ZT_j,|Q(T_j)|\sim\la_j,Q\in {\bf Q}[\mu_1,\mu_2]\}$. We fixed $\mu_j,\la_j$ in the rest of the argument and use $\ZT_j=\ZT_j[\la_j,\mu_1,\mu_2]$, ${\bf Q}={\bf Q}[\mu_1,\mu_2]$ for abbreviation.

Next, we introduce a collection of finitely overlapping $R^{1-\al}$-balls ${\bf B}=\{B\}$ in $B^{n+1}(0,2R)$ . For each ball $B\in {\bf B}$ and each tube $T_j\in\ZT_j$, consider the quantity
\begin{equation}
    \#\{Q\in{\bf Q}[\mu_1,\mu_2]:R^\de Q\cap T_j\not=\varnothing, Q\subset 2B\}.
\end{equation}
Let $B(T_j)$ be an $R^{1-\al}$ ball that maximizes the above quantity. We define a relation between $T_j\in\ZT_j$ and any $R^{1-\de}$ $B$ by: $T_j\sim B$ if and only if $B\subset3B(T_j)$. Thus, for any tube $T_j\in\ZT_j$,
\begin{equation}
\label{relation-Tao}
    \#\{B:T_j\sim B\}\lesssim 1.
\end{equation}
As a consequence, one has
\begin{equation}
\label{relation-consequence}
    \sum_B\#\{T_j:T_j\sim B,T_j\in\ZT_j\}=\sum_{T_j}\#\{B:T_j\sim B\}\lesssim |\ZT_j|.
\end{equation}
Using the relation defined above, we partition the left hand side of \eqref{kakeya-type} into four parts:
\begin{eqnarray}
\nonumber
    \Big|\sum_{T_1\in\ZT_1}\sum_{T_2\in\ZT_2}\phi_{T_1}\phi_{T_2}\Big|
    \!\!&=&\!\! \Big|\sum_B\Id_B\sum_{T_1\sim B}\sum_{T_2\sim B}\phi_{T_1}\phi_{T_2}+\sum_{T_1\sim B}\sum_{T_2\not\sim B}\phi_{T_1}\phi_{T_2}\\[1ex] \nonumber
    &&+\sum_{T_1\not\sim B}\sum_{T_2\sim B}\phi_{T_1}\phi_{T_2}+\sum_{T_1\not\sim B}\sum_{T_2\not\sim B}\phi_{T_1}\phi_{T_2}\Big|.
\end{eqnarray}
By triangle inequality, we only need to show
\begin{equation}
\label{related}
    \Big\|\sum_{Q\in{\bf Q}}\Id_Q\sum_B\Id_B\sum_{T_1\sim B}\sum_{T_2\sim B}\phi_{T_1}\phi_{T_2}\Big\|_{L^p(X)}\leq C_\e R^{(1-\al)\e}|\ZT_1|^{1/2}|\ZT_2|^{1/2}
\end{equation}
and for any subcollection $\ZT_2'\subset\ZT_2$
\begin{equation}
\label{unrelated}
    \Big\|\sum_{Q\in{\bf Q}}\Id_Q\sum_B\Id_B\sum_{T_1\not\sim B}\sum_{T_2\in\ZT_2'}\phi_{T_1}\phi_{T_2}\Big\|_{L^p(X)}\lesssim R^{C\al}|\ZT_1|^{1/2}|\ZT_2'|^{1/2},
\end{equation}
as if we pick $\al=\e^2$, the induction would close from $C_\e R^{(1-\al)\e}+R^{C\al}\lesssim R^\e$.

For \eqref{related}, we use the induction hypothesis \eqref{bilinear-Maximal-Schrodinger-esti} at scale $R^{1-\al}$, property (4) in Proposition \ref{wave-packet-paraboloid}, and \eqref{X-other-direction} so that
\begin{eqnarray}
\label{induction}
    &&\Big\|\sum_{Q\in{\bf Q}}\Id_Q\sum_B\Id_B\!\!\sum_{T_1\sim B}\sum_{T_2\sim B}\phi_{T_1}\phi_{T_2}\Big\|_{L^p(X)}\!\!\leq \sum_B \Big\|\sum_{T_1\sim B}\sum_{T_2\sim B}\phi_{T_1}\phi_{T_2}\Big\|_{L^p(B\cap X)}\\ \nonumber
    &&\leq C_\e R^{(1-\al)\e}\sum_B \big(\#\{T_1,T_1\sim B\})^{1/2}\big(\#\{T_2:T_2\sim B\}\big)^{1/2}.
\end{eqnarray}
Invoking Cauchy-Schwarz inequality, one can bounds \eqref{induction} as
\begin{equation}
\eqref{induction}\leq C_\e R^{(1-\al)\e} \big(\sum_B\#\{T_1,T_1\sim B\})^{1/2}\big(\sum_B\#\{T_2:T_2\sim B\}\big)^{1/2},
\end{equation}
which is 
\begin{equation}
    \lesssim C_\e R^{(1-\al)\e}|\ZT_1|^{1/2}|\ZT_2|^{1/2}.
\end{equation}
Thus, we finish the proof for \eqref{related}.

\medskip

It remains to prove \eqref{unrelated}. We need the following sharp $L^2$ estimates:
\begin{theorem}
\label{theoerm-du-zhang}
Let $f\in L^2(\ZR^n)$ with $\supp(\wh{f})\subset B^n(0,1)$. Then 
\begin{equation}
    \label{sharp-l2}
    \big\|\sup_{|t|<R} |\ce_{\Phi}f|\big\|_{L^2(B^n_R)}\leq C_{\e'} R^{\frac{n}{2(n+1)}+\e'}\|f\|_2.
\end{equation}
\end{theorem}
\noindent Theorem \ref{theoerm-du-zhang} was proved by Du and Zhang in \cite{Du-Zhang} for $\Phi(\xi)=|\xi|^2$. We will briefly explain why their argument works for general $\Phi$ in Remark \ref{remark1} at the end of this section.

Similar to \eqref{X-other-direction}, one has
\begin{eqnarray}
    && \Big\|\!\!\sum_{Q\in{\bf Q},B}\!\!\Id_Q\Id_B\!\!\sum_{T_1\not\sim B}\sum_{T_2\in\ZT_2'}\phi_{T_1}\phi_{T_2}\Big\|_{L^1(X)}\leq\Big\|\sum_{T_1\not\sim B}\phi_{T_1}\Big\|_{L^2(X)}\Big\|\sum_{T_2\in \ZT_2'}\phi_{T_2}\Big\|_{L^2(X)}\\ \nonumber
    &&\lesssim\Big\|\sup_{|t|<R} \big|\ce_{\Phi}\big(\sum_{T_1\not\sim B}f_{T_1}\big)\big|\Big\|_{L^2(B^n_R)}\Big\|\sup_{|t|<R} \big|\ce_{\Phi}\big(\sum_{T_2\in \ZT_2'}f_{T_2}\big)\big|\Big\|_{L^2(B^n_R)}\\ \nonumber
     &&+ \rapid(R)\|f_1\|_2^{1/2}\cdot\|f_2\|_2^{1/2}.
\end{eqnarray}
For simplicity, let us omit the error term $\rapid(R)\|f_1\|_2^{1/2}\cdot\|f_2\|_2^{1/2}$. Invoking \eqref{sharp-l2} we have
\begin{equation}
\label{l1}
    \Big\|\sum_{Q\in{\bf Q}}\Id_Q\sum_B\Id_B\sum_{T_1\not\sim B}\sum_{T_2\in\ZT_2'}\phi_{T_1}\phi_{T_2}\Big\|_{L^1(X)}\lesssim C_{\e'}R^{\frac{n}{n+1}+2\e'}|\ZT_1|^{1/2}|\ZT_2'|^{1/2}.
\end{equation}

On the other hand, Lee \cite{Lee-Bilinear} (Tao \cite{Tao-Bilinear} also provided a sketch of proof in the ending remark)  showed that 
\begin{equation}
\label{unrelated-l4}
    \Big\|\sum_{Q\in{\bf Q}}\Id_Q\sum_B\Id_B\sum_{T_1\not\sim B}\sum_{T_2\in\ZT_2'}\phi_{T_1}\phi_{T_2}\Big\|_{L^2(B^{n+1}_R)}\lesssim R^{C\al}R^{-\frac{n-1}{4}}|\ZT_1|^{1/2}|\ZT_2'|^{1/2}.
\end{equation}
Since $X\subset B^{n+1}_R$, we thus have
\begin{equation}
\label{l2}
    \Big\|\sum_{Q\in{\bf Q}}\Id_Q\sum_B\Id_B\sum_{T_1\not\sim B}\sum_{T_2\in\ZT_2'}\phi_{T_1}\phi_{T_2}\Big\|_{L^2(X)}\lesssim R^{C\al}R^{-\frac{n-1}{4}}|\ZT_1|^{1/2}|\ZT_2'|^{1/2}.
\end{equation}
Finally, we apply H\"older's inequality for \eqref{l1} and \eqref{l2} to get \eqref{unrelated}, and hence finish the proof for \eqref{bilinear-Maximal-Schrodinger-esti}. \qed

\begin{remark}
\label{remark1}
{\rm We briefly explain why Du-Zhang's argument remains valid for $\Phi$ satisfying $d(\Phi)=0$, $|\nabla\Phi|\lesssim1$, $|\det(\nabla^2\Phi)|\sim1$. Let us assume that $|\nabla\Phi|\leq C$ and every eigenvalue of $\nabla^2\Phi$ belongs to $[C^{-1},C]$. There are two reductions we would like to make:
\begin{itemize}
    \item By cutting the frequency unit ball $B^n(0,1)$ into smaller parts, we assume $|\nabla^3\Phi|\leq(100Cn)^{-100}$ without loss of generality. 
    \item After the first reduction, by rotation, we assume that $\nabla^2\Phi(0)$ is a diagonal matrix. 
\end{itemize}
As a result of the above two reductions, one can conclude via Taylor's theorem that every off-diagonal entry of $\nabla^2\Phi$ is bounded above by $(100Cn)^{-90}$.

We can run Du-Zhang's argument similarly in the "broad part", using multilinear refined Strichartz estimates \eqref{multilienar-refined-strichartz-esti}. After that, we need to be careful about the "narrow part": When $\Phi(\xi)=|\xi|^2$, the narrow part essentially lie in a lower dimension hyperplane (subspace). Since the intersection of this hyperplane with the graph $(\xi,|\xi|^2)$ is a lower dimension paraboloid, Bourgain-Demeter $l^2$ decoupling theorem is applicable. While for general $\Phi$, its narrow part lies in a thin neighborhood of the hypersurface determined by the equation
\begin{equation}
\label{hypersurface-gamma}
    m\cdot\nabla\Phi(\xi)+b=0
\end{equation}
for some vector $(m,b)\in S^{n+1}$ with $m\in\ZR^n,~|m|\geq 1/100$. In fact, let $\Ga=\Ga(\xi)$ be the hypersurface with codimension 2 in $\ZR^{n+1}$ obeying $m\cdot\nabla\Phi(\xi)+b=0$ and $\Phi(\xi)=0$. Then the narrow part falls in $N_{K^{-1}}(\Ga)$.

Since every eigenvalue of $\nabla^2\Phi$ belongs to $[C^{-1},C]$, the smooth map $\nabla^2\Phi:\ZR^n\to\ZR^n$ is morally an isometry. That is, for any unit vector $v\in\ZR^n$, one has $C^{-1}\leq|v\cdot\nabla^2\Phi|\leq C$. It in particular implies $|m\cdot\nabla^2\Phi(\xi)|\geq (100C)^{-1}$. After rotation, let us assume $m_k=0$ for $1\leq k\leq n-1$. By implicit function theorem, there is a smooth function $g(\xi_1,\ldots,\xi_{n-1})$ such that $m\cdot\nabla\Phi(\xi_1,\ldots,\xi_{n-1},g)+b=0$. One can use the first three derivative conditions on $\Phi$ to check $|\nabla g|\leq (100Cn)^{-80}$ and $|\nabla^2g|\leq (100Cn)^{-60}$. It implies that the hypersurface $\Ga$ has the parameterization $\Ga=\{(\bar\xi,g(\bar\xi),\Phi(\bar\xi,g(\xi)):\bar\xi\in B^{n-1}(0,1)\}$.

To prove a similar result as Lemma 3.3 in \cite{Du-Zhang}, we are going to use a bootstrapping argument. When $\Phi$ is the paraboloid $|\xi|^2$, one can use Bourgain-Demeter's $l^2$ decoupling theorem to decouple the narrow part into $K^{-1}$ caps directly. This is what Du-Zhang did in \cite{Du-Zhang}. In our case, it is hard to do that since the neighborhood $N_{K^{-1}}(\Ga)$ is bent in the $n$-th direction. While we can still use the decoupling theorem at a much larger scale $c:=(100Cn)^{-40}$, since 
the neighbourhood $N_{K^{-1}}(\Ga)$ is contained in a thin rectangular stripe $\{(\xi,\Phi(\xi)):|\xi_n|\leq c\}$, which has width $c$. Thus, we use the decoupling theorem to decouple $N_{K^{-1}}(\Ga)$ into $c^{-{(n-1)}}$ pieces, each of which lies in a $c$-cap. Note that in each $c$-cap, the neighbourhood $N_{K^{-1}}(\Ga)$ is contained in an even smaller region---a thin rectangular stripe of width $c^2$. Thus, after parabolic rescalings, we can again use the decoupling theorem to decouple each $c$-cap into $c^{-(n-1)}$ many $c^2$-caps. We repeat our argument until we have a collection of $c^{2^\ka}$-caps, where
$\ka:=[\log_2\log_{c^{-1}} K]\sim \log\log K$. At this point, the original neighbourhood $N_{K^{-1}}(\Ga)$ is decoupled into $K^{-1}$-caps as desired.

We also need to track our loss at each step when using the decoupling theorem. Indeed, when using the decoupling theorem to decouple $c^{2^j}$-caps into $c^{2^{j+1}}$-caps, the loss can be as small as $C_\e c^{-2^j\e}$. Hence the total loss is bounded above by $\prod_{j=1}^\ka C_\e c^{-2^j\e}\lesssim c^{-2^\ka\e}\lesssim K^\e$, which is an acceptable one.

}
\end{remark}

\section{Appendix: Proof of wave-packet decomposition }
In this section we present a proof for Proposition \ref{wave-packet-paraboloid}. The proof based on the framework in \cite{Guth-restriction-R3} Proposition 2.6. Let $\psi(\xi)$ be a smooth function that equals to $1$ in the unit ball $B^n(0,1)$, and is supported in a bigger ball $B^n(0,2)$. Define $\psi_q(\xi):=\psi\big(R^{1/2}(\xi-c(q))\big)$ so that $\wh{f}_q=\psi_q\wh{f}_q$. Consider the partial Fourier series $S_N\wh{f_q}$ for $\wh{f_q}$ expanding in a $2R^{-1/2}$-cube $2q$:
\begin{equation}
    S_N\wh{f_q}(\xi)\sim\sum_{m\in\ZZ^n,|m|\leq N}a_me^{i\pi R^{1/2}m\cdot\xi},
\end{equation}
where 
\begin{equation}
    a_m:=R^{n/2}\int_{2q} e^{-i\pi R^{1/2}m\cdot\xi}\wh{f_q}(\xi)\psi_q(\xi)d\xi.
\end{equation}

For $T\in\ZT_q$, let $P_T(x):\ZR^{n+1}\to\ZR^n$ be the projection to the subspace whose normal vector coincides to the direction of $T$. If $P_T(c(T))=R^{1/2}m$, we let $m=m_T$ and define
\begin{equation}
    f_T(\xi)=\int a_m e^{ix\cdot\xi}e^{i\pi R^{1/2}m\cdot\xi}\psi_q(\xi)d\xi
\end{equation}
so that $\wh{f_T}(\xi)=a_me^{i\pi R^{1/2}m\cdot\xi}\psi_q(\xi)$. Clearly, Property (1) is true.

Next we take a look on $\ce_{\Phi} f_T$. Plug in the definition of $f_T$ to have
\begin{eqnarray}
\nonumber
    \ce_{\Phi} f_T(x,t) \!\!\!&=&\!\!\! \int_{\ZR^n}e^{ix\cdot\xi}e^{it\Phi(\xi)}a_me^{i\pi R^{1/2}m\cdot\xi}\psi_q(\xi)d\xi\\ \nonumber
    &=&\!\!\!  \int_{\ZR^n}e^{ix\cdot\xi}e^{it\Phi(\xi)}a_me^{i\pi R^{1/2}m\cdot\xi}\psi(R^{1/2}(\xi-c(q)))d\xi.
\end{eqnarray}
After a change of variable, we use Taylor's expansion so that
\begin{equation}
\label{one-pack}
    |\ce_{\Phi} f_T(x,t)|= \Big|a_m\int_{\ZR^n}e^{i(x+\pi R^{1/2}m+\nabla\Phi(c(q)))\cdot\xi}e^{iO(1)}\psi(R^{1/2}\xi)d\xi\Big|.
\end{equation}
When $0<|t|<R$ and $|x+\pi m+\nabla\Phi(c(q))|\gtrsim R^{1/2}R^\de$, the integrand in \eqref{one-pack} admits fast decay. Thus $| \ce_{\Phi} f_T|\lesssim \rapid(R)\|f\|_2$, as $|a_m|\leq R^{n/4}\|f\|_2$. This gives the proof of Property (2).

We now prove Property (3). Recall that the partial sum $S_N\wh{f_q}$ converges to $\wh{f_q}$ in $L^2$. Thus, there is a positive number $N_q>0$ such that
\begin{equation}
\label{l2-convergence}
    \|S_{N_q}\wh{f_q}-\wh{f_q}\|_{L^2(2q)}\leq \rapid(R)\|f\|_2.
\end{equation}
Let $\bar\ZT_q$ be the collection of $T\in\ZT_q$ such that $|P_T(c(T))|\leq N_q$. Then
\begin{eqnarray}
\label{partial-sum-esti}
    && \Big|\ce_{\Phi}f_q-\sum_{T\in\ZT_q}\ce_{\Phi}f_T\Big| = \Big|\int_{\ZR^n} e^{ix\cdot\xi}e^{it\Phi(\xi)}(\wh{f_q}-\sum_{T\in\ZT_q} \wh{f_T})d\xi\Big|\\ \nonumber
    && \leq \Big|\int_{\ZR^n} e^{ix\cdot\xi}e^{it\Phi(\xi)}(\wh{f_q}-\sum_{T\in\bar\ZT_q} \wh{f_T})d\xi\Big|+\Big|\int_{\ZR^n} e^{ix\cdot\xi}e^{it\Phi(\xi)}\Big(\sum_{T\in(\bar\ZT_q\setminus\ZT_q)} \wh{f_T}\Big)d\xi\Big|.
\end{eqnarray}
From \eqref{l2-convergence} and the fact $\wh{f}_q=\psi_q\wh{f}_q$, the first part of \eqref{partial-sum-esti} is bounded above by $\rapid(R)\|f\|_2$. Since $(x,t)\in B^{n+1}_R$, for each $T\in\bar\ZT_q\setminus\ZT_q$, we set $m=P_T(C(T))$ and use the standard (non) stationary phase method to get
\begin{equation}
    \Big|\int_{\ZR^n} e^{ix\cdot\xi}e^{it\Phi(\xi)} \wh{f_T}(\xi)d\xi\Big|\lesssim \frac{a_m}{(1+|x|+|m|)^{-M}}\lesssim |m|^{-2n}\rapid(R)\|f\|_2
\end{equation}
for some $M$ large enough. Summing up all the $T\in\bar\ZT_q\setminus\ZT_q$ we have the second part of \eqref{partial-sum-esti} is bounded by $\rapid(R)\|f\|_2$. Thus, \eqref{partial-sum-esti} is bounded by $\rapid(R)\|f\|_2$ and we finish the proof for Property (3) by summing up all the $q$.

Property (4), the first part of Property (5) and Property (6) follow directly from Plancherel. For the second part of Property (5), by Plancherel we have
\begin{equation}
    \langle \ce_{\Phi} f_T,\ce_{\Phi} f_{T'}\rangle=\int \langle\wh{f_T},\wh{f}_{T'}\rangle dt\leq R|\langle\wh{f_T},\wh{f}_{T'}\rangle|,
\end{equation}
which is further bounded by
\begin{equation}
    CR^{n/2}|a_{m_T}a_{m_{T'}}\wh\psi(m_T-m_{T'})|\lesssim\rapid(R),
\end{equation}
as $|m_T-m_{T'}|\sim R^{-1/2}${\rm dist($T,T'$)}$\gtrsim R^{\de}$. \qed

\section{Appendix: An epsilon removal lemma}
In this section we prove the following lemma
\begin{lemma}
\label{epsilon-removal}
Suppose $p>2$, $\e>0$ and \eqref{Maximal-Schrodinger-local-esti}. Then 
\begin{equation}
    \big\|\sup_{t\in\ZR }|\ce_{\Phi}f|\big\|_{L^{p_0}(\ZR^n)}\lesssim \|f\|_2
\end{equation}
for $f\in L^2(\ZR^n)$ and
\begin{equation}
\label{removal-range}
    \frac{1}{p_0}<\frac{1}{p}-\frac{C}{\log \frac{1}{\e}}.
\end{equation}
\end{lemma}
In particular, Lemma \ref{Maximal-Schrodinger-local} implies Theorem \ref{Maximal-Schrodinger} by letting $\e\to0$. Our proof for Lemma \ref{epsilon-removal} is similar to the argument in \cite{Bourgain-Guth-Oscillatory}. See also \cite{Tao-BR-Restriction}. 

For some technical issues, we assume $\wh{f}$ is supported in $B^n_{1/4}$ rather then the unit ball in Theorem \ref{Maximal-Schrodinger}. Since $\wh{f}\subset B^n_{1/4}$, $|\ce_{\Phi}f|$ is essentially constant in every 1-ball in $\ZR^{n+1}$. Inspired by this observation, we have the following lemma:
\begin{lemma}
\label{weight}
There exists a function $\vp(x,t):\ZR^{n+1}\to\ZC$ such that
\begin{equation}
\label{weight-repre}
    \big\|\sup_{t\in\ZR}|\ce_{\Phi}f|\big\|_{L^p(\ZR^n)}\sim\big\|(\ce_{\Phi}f)\vp\big\|_{L^p(\ZR^{n+1})},
\end{equation}
where the function $\vp$ satisfies the following properties:
\begin{enumerate}
    \item  $|\vp|\lesssim1$, ${\rm supp}(\wh\vp)\subset B^{n+1}_{1/4}$ .
    \item Uniformly for any $x\in \ZR^n$, $\|\vp(x,\cdot)\|_{L^1_t}=O(1)$.
    \item For any small factor $\be>0$, there exists $R^\be$ many horizontally sparse sets $\{X_j\}$ (See Definition \ref{set-X}) such that $\vp(x,t)=\rapid(R)$ when $(x,t)\in\ZR^{n+1}\setminus(\cup_j X_j)$. 
\end{enumerate}
We remark that the implicit constant in $\rapid(R)$ depends on $\be$
\end{lemma}
Roughly speaking, the weight $\vp$ is an averaging method to help us realize the linearization $\sup_t|\ce_\Phi f(x,t)|=|\ce_\Phi f(x,t(x))|$. We remark that $\vp$ is essentially supported in a set satisfying Property \ref{set-X}.

\begin{proof}
Let $\{U\}$ be the lattice $1$-cubes in $\ZR^n$ and $\{I\}$ be the lattice $1$-cubes in $\ZR$. Let $\psi_U(x)$, $\psi_I(t)$ be two smooth functions on such that ${\rm supp}(\wh\psi_U)\subset B^n_{1/8}$, ${\rm supp}(\wh\psi_I)\subset[-1/8,1/8]$, $|\psi_U(x)|\sim1$ for $x\in U$, $|\psi_I(t)|\sim1$ for $t\in I$ and $\psi_U$, $\psi_I$ admit fast decay outside $U$, $I$, respectively. Thus, by Hausdorff-Young inequality,
\begin{equation}
\label{L1-Linfty}
    \int_U\sup_{t\in\ZR}|\ce_{\Phi} f|^p\lesssim \sup_{I}\sup_{x\in\ZR^n}\sup_{t\in I}\big|\ce_{\Phi} f(\psi_U\psi_I)^3\big|^p\leq\sup_I\big\|\big(\ce_{\Phi} f(\psi_U\psi_I)^3\big)^\wedge\big\|_1^p.
\end{equation}
From the constructions of $\psi_U$ and $\psi_I$, we see the Fourier transform of $(\psi_U\psi_I)^3$ is compactly supported. Combining the fact that $\wh{f}$ is compactly supported,  we have that $\big(\ce_{\Phi} f(\psi_U\psi_I)^3\big)^\wedge$ is compactly supported. Hence
\begin{equation}
\label{L2-L2}
    \sup_I\big\|\big(\ce_{\Phi} f(\psi_U\psi_I)^3\big)^\wedge\big\|_1^p\lesssim\sup_I\big\|\big(\ce_{\Phi} f(\psi_U\psi_I)^3\big)^\wedge\big\|_2^p.
\end{equation}
Invoking Plancherel and H\"older's inequality, we obtain that for $3<p<4$,
\begin{eqnarray}
\label{Holder}
    &&\sup_I\big\|\big(\ce_{\Phi} f(\psi_U\psi_I)^3\big)^\wedge\big\|_2^p=\sup_I\big\|\ce_{\Phi} f(\psi_U\psi_I)^3\big\|_2^p\\ \nonumber
    &&\leq\sup_I\Big(\int_{\ZR^{n+1}}|\ce_{\Phi} f|^p|\psi_U\psi_I|^{2p}\Big)\Big(\int_{\ZR^{n+1}}|\psi_U\psi_I|^{2p/(p-2)}\Big)^{(p-2)/2}\\ \nonumber
    &&\lesssim \sup_I\int_{\ZR^{n+1}}|\ce_{\Phi} f|^p|\psi_U\psi_I|^{2p}.
\end{eqnarray}
We pick one $I_U$ such that 
\begin{equation}
\label{discretize}
    \sup_I\int_{\ZR^{n+1}}|\ce_{\Phi} f|^p|\psi_U\psi_I|^{2p}\leq 2\int_{\ZR^{n+1}}|\ce_{\Phi} f|^p|\psi_U\psi_{I_U}|^{2p}
\end{equation}
and define
\begin{equation}
    \vp:=\sum_U|\psi_U\psi_{I_U}|^2.
\end{equation}
Note that $\psi_U(x)\psi_{I_U}(t)=\rapid(R)$ when $\dist((x,t),U\times I_U)\geq R^\be$ for any $\be>0$. One can check directly that the function $\vp$ satisfies all properties mentioned in the statement of the Lemma. Combining \eqref{L1-Linfty}, \eqref{L2-L2}, \eqref{Holder}, \eqref{discretize} and summing all the $U\subset\ZR^n$, we have
\begin{equation}
    \int_{\ZR^n}\sup_{t\in\ZR}|\ce_{\Phi} f|^p\lesssim\sum_U|\int_{\ZR^{n+1}}|\ce_{\Phi} f|^p|\psi_U\psi_{I_U}|^{2p}\leq\int_{\ZR^{n+1}}|\ce_{\Phi}\vp|^p.
\end{equation}
Take $p$-th root to both sides we get one direction for \eqref{weight-repre}. For the other direction of the estimate above, one just need to use property (2) of the function $\vp$.
\end{proof}

Applying Lemma \ref{weight}, it suffices to show the dual estimate
\begin{equation}
\label{dual-global}
    \|\ce^\ast (g\vp)\|_{L^2(B^{n}_{1/4})}\lesssim\|g\|_{p_0'}
\end{equation}
where 
\begin{equation}
    \ce^\ast g(\xi)=\int_{\ZR^{n+1}} e^{-ix\cdot\xi}e^{-it\Phi(\xi)}g(x,t)dxdt.
\end{equation}

Let us follow a similar argument as in the proof of Lemma \ref{weight}. Assuming \eqref{Maximal-Schrodinger-local-esti}, we have that for any $R$-ball $V\subset\ZR^{n+1}$, 
\begin{equation}
    \big\|(\ce_{\Phi}f)\vp\big\|_{L^p(V)}\leq C_\e R^\e\|f\|_2
\end{equation}
and the dual estimate
\begin{equation}
\label{dual-local}
    \|\ce^\ast (g\Id_{V}\vp)\|_{L^2(B^{n}_{1/4})}\leq C_\e R^\e\|g\|_{p'}.
\end{equation}

Let $\phi(x)$ be a smooth function in $\ZR^n$ that $\wh\phi(\xi)=1$ when $\xi\in B^n_{1/4}$ and ${\rm supp}(\wh\phi)\subset B^n_{1/2}$. The classic result of restriction theorem (see \cite{Hormander} Section 7, \cite{Stein-Oscillatory} Proposition 6) tells us that $\ce_{\Phi}\phi(x,t)$ is bounded above by the decay function $C(1+|x|+|t|)^{-n/2}$. Therefore, in order to make full use of the local estimate \eqref{dual-local}, we are motivated to consider a sparse collection of $R$-ball in $\ZR^{n+1}$.

Let $\{V_j\}_{j=1}^N$ be a collection of $R$-balls in $\ZR^{n+1}$ such that for any $j\not=j',j,j'\in\{1,2,\ldots,N\}$, ${\rm dist}\big(c(V_j),c(V_{j'})\big)\gtrsim R^{2C}N^{2C}$. Here $C$ is a large absolute constant that will be determined later. For a measurable function $g:\ZR^{n+1}\to\ZC$, we let $G=\sum_{j}g\Id_{V_j}$. Consider $\|\ce^\ast (G\vp)\|_2$:
\begin{eqnarray}
\label{expand-esti}
    &&\big\|\ce^\ast G\vp\big\|_{L^2(B^{n}_{1/4})}^2\leq\int_{\ZR^n} \Big|\int_{\ZR^{n+1}} e^{-ix\cdot\xi}e^{-it\Phi(\xi)}(G\vp)(x,t)dxdt\Big|^2\phi(\xi)d\xi\\ \nonumber
    &&= \int_{\ZR^{n+1}}\!\int_{\ZR^{n+1}}\!\! \Big(\int_{\ZR^n} e^{i(y-c)\cdot\xi}e^{i(s-t)\Phi(\xi)}\phi(\xi)d\xi\Big)(G\vp)(x,t)\overline{(G\vp)}(y,s)dxdtdyds.
\end{eqnarray}
Extract the kernel
\begin{equation}
    \label{kernel}
    K(x-y,t-s)=\int_{\ZR^n} e^{i(y-x)\cdot\xi}e^{i(s-t)\Phi(\xi)}\phi(\xi)d\xi
\end{equation}
and expand $G=\sum_{V_j}g\Id_{V_j}$ in the right hand side of \eqref{expand-esti} so that
\begin{eqnarray}
\label{TT*}
    \big\|\ce^\ast (G\vp)\big\|_{L^2(B^{n}_{1/4})}^2 \!\!\!\!\!&\leq&\!\!\!\! \int_{\ZR^{2n+2}}\!\!\! K(x\!-\!y,t\!-\!s)\sum_{j,j'}(g\Id_{V_j}\vp)(x,t)(\overline{g\Id_{V_{j'}}\vp})(y,s)\\ \nonumber
    &=&\!\!\!\! \sum_{j,j'}\int_{\ZR^{2n+2}} \!\!\!K(x\!-\!y,t\!-\!s)(g\Id_{V_j}\vp)(\overline{g\Id_{V_{'}}\vp}).
\end{eqnarray}

For $j=j'$, we use \eqref{dual-local} to have
\begin{equation}
\label{diagonal}
    \int K(x-y,t-s)(g\Id_{V_j}\vp)(\overline{g\Id_{V_{j'}}\vp})\leq\big\|\ce^\ast (g\Id_{V_j}\vp)\big\|_{L^2(B^n_{1/2})}^2\leq C_\e R^\e\|g\Id_{V_j}\|_{p'}^2.
\end{equation}
For $j\not=j'$, we apply H\"older inequality to get
\begin{equation}
\label{off-diagonal1}
   \int_{\ZR^{2n+2}} K\cdot(g\Id_{V_j}\vp)(\overline{g\Id_{V_{j'}}\vp}) \leq \|g\Id_{V_j}\|_{p'}\Big\|\int \wt{K}(x,t,y,s)(\overline{g\Id_{V_j}})(y,s)dyds\Big\|_p
\end{equation}
where
\begin{equation}
    \wt{K}(x,t,y,s)=\vp(x,t)\Id_{V_j}(x,t)K(x\!-\!y,t\!-\!s)\overline{\vp(y,s)}\Id_{V_{j'}}(y,s).
\end{equation}
Notice that $K(x-y,t-s)=\ce_{\Phi}\phi(x-y,t-s)$ and $|\ce_{\Phi}\phi(x,t)|\lesssim(1+|x|+|t|)^{-n/2}$. Thus, together with $|\vp|\lesssim1$ and generalized Young's inequality we have
\begin{equation}
\label{off-diagonal2}
    \Big\|\int \wt{K}\cdot(\overline{g\Id_{V_j}})dyds\Big\|_p\lesssim\|K\Id_{\{V_j-V_{j'}\}}\|_{p/2}\|(g\Id_{V_j})\|_{p'}\lesssim (NR)^{-C}\|(g\Id_{V_j})\|_{p'}.
\end{equation}
Combining \eqref{TT*}, \eqref{diagonal}, \eqref{off-diagonal1} and \eqref{off-diagonal2}, we get
\begin{eqnarray}
\label{sparse}
     \big\|\ce^\ast G\big\|_{L^2(B_{1/4}^n)}^2 \!\!\!\!&\lesssim& \!\!\! C_\e R^\e\sum_{V_j}\|g\Id_{V_j}\|_{p'}^2+\sum_{j\not=j'}N^{-C}R^{-C}\|g\Id_{V_j}\|_{p'}\|g\Id_{V_{j'}}\|_{p'}\\ \nonumber
     &\lesssim& \!\!\! C_\e R^\e\Big(\big\|\sum_{V_j}g\Id_{V_j}\big\|_{p'}^{p'}\Big)^{2/p'}\leq C_\e' R^\e\|G\|_{p'}^2.
\end{eqnarray}

We say a set $E=\cup_{j=1}^N B^{n+1}(x_j,r)$ is \emph{ sparse} if ${\rm dist}(x_j,x_{j'})\gtrsim r^{2C}N^{2C}$. Therefore, what we have proved above is the following lemma:
\begin{lemma}
Let $\{V_j\}_{j=1}^N$ be a collection of sparse $R$-balls in $\ZR^{n+1}$ and $G=\sum_{j}g\Id_{V_j}$, then
\begin{equation}
\label{sparse-local}
    \big\|\ce^\ast G\big\|_{L^2(B_{1/4}^n)}\leq C_\e R^\e\|G\|_{p'}^2.
\end{equation}
\end{lemma}

We will use \eqref{sparse-local} to prove \eqref{dual-global}. Since $\wh{\vp}$ is supported on $B^{n+1}_{1/2}$, without loss of generality, we assume $\wh{g}$ is supported in $B^{n+1}_{7/8}$. Let $\wh\chi$ is a smooth function supported on $B^n_1$ and $\wh{\chi}(\xi)=1$ on $B^{n+1}_{7/8}$, then $g=\chi\ast g=(\chi\ast\chi)\ast g$. We also assume $g$ is a Schwartz function so that its Fourier series converges. We expand $\wh{g}$ on $B^{n+1}_1$ to have
\begin{equation}
    \wh{g}(\xi,\tau)=\wh\chi^2(\xi,\tau)\wh{g}(\xi,\tau)\sim\sum_{(m,n)\in\ZZ^{n+1}}a_{m,n}e^{i\pi m\cdot\xi}e^{i\pi n\tau}\wh{\chi}(\xi,\tau)
\end{equation}
with
\begin{equation}
    a_{m,n}=\int_{\ZR^{n+1}} e^{-i\pi m\cdot\xi}e^{-i\pi n\tau}\wh{g}(\xi,\tau)\wh{\chi}(\xi,\tau)d\xi d\tau\sim\chi\ast g(-m,-n).
\end{equation}

Next, we assume $\sum_{m,n}|a_{m,n}|^{p_0}=1$, and sort $|a_{m,n}|$ dyadically. Let 
\begin{equation}
    g_k(x,t)=\sum_{(m,n)\in\La_k}a_{m,n}\chi(x+m,t+n)
\end{equation}
where $\La_k=\{(m,n):|a_{m,n}|\sim2^{-k}\}$ so that $|\La_k|\lesssim 2^{kp_0'}$. Since
\begin{enumerate}
    \item $\|g_k\|_1\leq \sum|a_{m,n}|\sim 2^{-k}|\La_k|$,
    \item $\|g_k\|_2\sim(\sum|a_{m,n}|^2)^{1/2}\sim 2^{-k}|\La_k|^{1/2}$,
    \item $\|g_k\|_\infty\leq \sup_{m,n}|a_{m,n}|\lesssim 2^{-k}$,
\end{enumerate}
applying H\"older twice, for $1<p'<2$, we have
\begin{equation}
    \|g_k\|_{p'}\sim 2^{-k}|\La_k|^{1/p'}.
\end{equation}

We need a covering lemma in \cite{Tao-BR-Restriction}.
\begin{lemma}
\label{covering-lemma}
Suppose $E$ is a union of $1$-cubes. Then there exist $O(N|E|^{1/N})$ collections of sparse set that cover $E$, such that the balls in each sparse set have radius at most $O(|E|^{C^N})$.
\end{lemma}

Now we are in a position to prove \eqref{dual-global}. Since $g=\sum_{k\geq0}g_k$, by triangle inequality
\begin{equation}
\label{triangle-ineq}
    \|\ce^\ast (g\vp)\|_{L^2(B^n_{1/4})}\leq\sum_{k\geq0}\|\ce^\ast (g_k\vp)\|_{L^2(B^n_{1/4})}.
\end{equation}
For each $k\geq0$, since the function $\chi$ admits fast decay, $g_k$ is essentially supported on a set $E$, where $E$ is a union of $1$-cubes, and $|E|\lesssim|\La_k|\lesssim 2^{kp_0'}$. We apply Lemma \ref{covering-lemma} so that
we can obtain a collection of sparse sets ${\bf E}=\{E_j\}$ such that $E\subset\cup E_j$ and $|{\bf E}|\lesssim N|E|^{1/N}$. By triangle inequality and \eqref{sparse-local},
\begin{equation}
\nonumber
    \big\|\ce^\ast (g_k\vp)\big\|_{L^2(B^n_{1/4})}\lesssim\sum_j\big\|\ce^\ast (2^{-k}\Id_{E_j\cap E}\vp)\big\|_{L^2(B^n_{1/4})}\leq 2^{-k} C_\e |E|^{\e C^N}N|E|^{\frac{1}{N}}|E|^{\frac{1}{p'}},
\end{equation}
which is further bounded by
\begin{equation}
\label{bound}
    2^{-k}C_\e\frac{\log(1/\e)}{C}|E|^{\e^{1-\frac{\log C}{C}}+\frac{C}{\log(1/\e)}+\frac{1}{p'}},
\end{equation}
if we let $N=\log(1/\e)/C$. Since $|E|\lesssim 2^{kp_0'}$, plug it in we have
\begin{equation}
    \eqref{bound}\lesssim C_\e\frac{\log(1/\e)}{C}2^{k(p_0'\e^{1-\frac{\log C}{C}}+\frac{p_0'C}{\log(1/\e)}+\frac{p_0'}{p'}-1)}.
\end{equation}
Thus, if first choose $C$ big enough then $\e$ small enough and let $p_0$ be in \eqref{removal-range}, we get that for some small positive number $\e'$,
\begin{equation}
    \eqref{bound}\lesssim C_{\e'}2^{-k\e'}.
\end{equation}
Summing up all the $k$ we therefore can conclude
\begin{equation}
\label{discrete-esti}
    \|\ce^\ast (g\vp)\|_{L^2(B^n_{1/4})}\lesssim_{p_0}1=\Big(\sum_{m,n}|a_{m,n}|^{p_0}\Big)^{1/p_0}.
\end{equation}
Noticing that \eqref{discrete-esti} is also true with $a_{m,n}$ replaced by $a_{m+x,n+t}$ where $(x,t)\in B^{n+1}_1$ and
\begin{equation}
\nonumber
    a_{m+x,n+t}=\int_{\ZR^{n+1}} e^{-i(m+x)\cdot\xi}e^{-i(n+t)\tau}\wh{g}(\xi,\tau)\wh{\chi}(\xi,\tau)d\xi d\tau\sim \chi\ast g(-m\!-\!x,-n\!-\!t),
\end{equation}
we can average over the translated $\ZZ^{n+1}$-lattices so that
\begin{eqnarray}
&& \|\ce^\ast (g\vp)\|_{L^2(B^n_{1/4})}=\int_{B^{n+1}_1}\|\ce^\ast (g\vp)\|_{L^2(B^n_{1/4})}dxdt \\ \nonumber
&& \lesssim \int_{B^{n+1}_1}\Big(\sum_{m,n}|a_{m+x,n+t}|^{p_0}\Big)^{1/p_0'}dxdt,
\end{eqnarray}
which is 
\begin{eqnarray}
    \sim\int_{B^{n+1}_1}\Big(\sum_{m,n}|\chi\ast g(-m\!-\!x,-n\!-\!t)|^{p_0}\Big)^{1/p_0'}dxdt
\end{eqnarray}
and is further bounded by
\begin{equation}
\label{averaged-eq}
    \Big(\int_{B^{n+1}_1}\sum_{m,n}|\int_{\ZR^{n+1}}\chi(-x\!-\!m\!-\!y,-t\!-\!n\!-s)g(y,s)dyds|^{p_0'}dxdt\Big)^{1/p_0'}.
\end{equation}

Finally, since $\chi\ast g=g$ and \eqref{averaged-eq} is nothing but $\|\chi\ast g\|_{p_0'}$, we get \eqref{dual-global} and hence finish the proof of Lemma \ref{epsilon-removal}. \qed

\end{document}